\documentclass[final,1p,times,number]{elsarticle}

\usepackage{amsmath}
\usepackage{amssymb}
\usepackage{amsfonts}
\usepackage{graphicx}
\usepackage{epstopdf}
\usepackage{color}
\usepackage{bm}
\usepackage{multirow}
\usepackage{natbib}
\usepackage{hyperref}
\usepackage{cleveref}

\usepackage{color}
\usepackage{float}
\usepackage{caption}

\usepackage[caption=false]{subfig}

\usepackage{ifpdf}
\ifpdf%
\usepackage{pdflscape}
\else
\usepackage{lscape}
\fi

\hypersetup{
	colorlinks,
	linkcolor={red!50!black},
	citecolor={blue!50!black},
	urlcolor={blue!80!black}
}

\def\a{\mbox{\boldmath $a$}}

\def\e{\mbox{\boldmath $e$}}
\def\f{\mbox{\boldmath $f$}}

\def\g{\mbox{\boldmath $g$}}
\def\h{\mbox{\boldmath $h$}}
\def\m{\mbox{\boldmath $m$}}

\def\x{\mbox{\boldmath $x$}}
\def\y{\mbox{\boldmath $y$}}

\def\0{\mbox{\boldmath $0$}}

\biboptions{sort&compress}

\begin{document}

\begin{frontmatter}
% \title{An Efficient Energy Stable Structure Preserving Method for Micromagnetics Simulations}

\title{An Efficient Energy Stable Structure Preserving Method for The Landau-Lifshitz Equation}

\author[XJTLU]{Changjian Xie\corref{cor}}
\cortext[cor]{Corresponding author.} 
\ead{Changjian.Xie@xjtlu.edu.cn}

\author[XJTLU]{Yingxi Miao}
\ead{Yingxi.Miao22@student.xjtlu.edu.cn}

\author[XJTLU]{Haocheng Yang}
\ead{Haocheng.Yang22@student.xjtlu.edu.cn}

\address[XJTLU]{School of Mathematics and Physics, Xi'an-Jiaotong-Liverpool University, Re'ai Rd. 111, Suzhou, 215123, Jiangsu, China.}

\begin{abstract}
One of the main difficulties in micromagnetics simulation is the norm preserving constraints $\|\m\|=1$ at the continuous or the discrete level. Another difficulty is the stability with the time step constraint. Using standard explicit integrators leads
to a physical time step of sub-pico seconds, which is often two orders of magnitude
smaller than the fastest physical time scales. Direct implicit integrators require solving complicated, coupled systems. Another major difficulty with the projection method in this field is the lack of rigorous theoretical guarantees regarding its stability of the projection step. In this paper, we introduce an first order method. Such a
method is structure preserving based on a combination of a Gauss–Seidel iteration, a double diffusion iteration and a Crank-Nicolson iteration to preserve the norm constraints. 
% This method allows us to carry out fully resolved calculations
% for the switching of the magnetization in micron-sized elements. 
\end{abstract}


\end{frontmatter}

\section{Introduction}

Ferromagnetic materials find widespread application in data storage technologies, largely due to the bistable nature of their intrinsic magnetic ordering, commonly referred to as magnetization. The temporal evolution of magnetization in these systems is principally dictated by the Landau-Lifshitz-Gilbert (LLG) equation \cite{Landau1935On,Gilbert:1955}. This equation incorporates two fundamental terms governing magnetization dynamics: a gyromagnetic term, responsible for energy conservation, and a damping term, which models energy dissipation. 
% The damping term is critically important in magnetic devices, profoundly influencing both their energy consumption and operational speed. Recent experimental work on magnetic-semiconductor heterostructures \cite{Zhang2020ExtremelyLM} has confirmed that the Gilbert damping constant can be tuned. At the microscopic scale, damping arises from various physical processes, including electron scattering, itinerant electron relaxation \cite{Heinrich1967TheIO}, and phonon-magnon coupling \cite{Suhl1998TheoryOT, Nan2020ElectricfieldCO}, which can be quantitatively assessed through electronic structure calculations \cite{TangXia2017}. From an applied standpoint, modulating the damping parameter enables optimization of magnetodynamic properties in ferromagnetic materials—for instance, by lowering the switching current and enhancing the writing speed of magnetic memory devices \cite{Wei2012MicromagneticsAR}. Although many experimental studies have focused on systems with small damping parameters \cite{Budhathoki2020LowGD,Lattery2018LowGD,Weber2019GilbertDO}, substantial damping effects have also been reported in works like \cite{GilbertKelly1955, Tanaka2014MicrowaveAssistedMR}. Notably, Tanaka et al. \cite{Tanaka2014MicrowaveAssistedMR} found that larger damping constants correlate with reduced magnetization switching times. Gilbert and Kelly \cite{GilbertKelly1955} also reported observing very large damping parameters (around 9) in their investigations. 
The LLG equation constitutes a vector-valued, nonlinear system characterized by a point-wise constant magnitude constraint on the magnetization vector. Significant research efforts have been dedicated to devising efficient and numerically stable methods for micromagnetic simulations, as summarized in review articles such as \cite{kruzik2006recent,cimrak2007survey}. 

To handle the point-wise constant magnitude constraint on the magnetization vector, the methods are mainly divided into the following categories:
\begin{itemize}
    \item Projection methods: By projecting the numerical solution onto the unit sphere, the constraint of constant modulus is maintained, with representative work such as GSPM \cite{wang2001gauss}.Although this method is simple for engineering applications, spherical projection is a nonlinear operation, which brings certain challenges to numerical analysis. Review paper refer to \cite{Shen2006} and reference therein.
    \item Lagrange multiplier methods: such a method constructs efficient, robust higher-order schemes with different Lagrange multiplier via a predictor-corrector method, and further develops length-preserving and energy-dissipative schemes at the cost of solving one nonlinear algebraic equation \cite{qing_2023}. Such a method need to use a construction for extra Lagrange multiplier. Another method is generalized scalar auxiliary variable (GSAV) approach to preserve the constraints, see \cite{xiaoli_2024}.
    \item Fully implicit methods: such a method have great stability and accuracy. However, the nonlinear system is resolved by \cite{JEONG2010613}, \cite{An2016OptimalEE}.
\end{itemize}

% Semi-implicit schemes have attracted considerable interest among existing numerical strategies due to their ability to maintain good numerical stability without requiring complex nonlinear solvers \cite{alouges2006convergence, gao2014optimal, Xie2018}. For instance, our research group previously proposed a second-order backward differentiation formula (BDF2) scheme based on one-sided interpolation \cite{Xie2018}. This method requires solving a three-dimensional linear system with non-constant coefficients at each time step. A rigorous theoretical analysis establishing the second-order convergence of this BDF2 method was provided in \cite{jingrun2019analysis}. Alternatively, Alouges et al. \cite{alouges2006convergence}  introduced a linearly implicit method utilizing the tangent space to preserve the magnetization constraint, although this approach achieves only first-order temporal accuracy. More recently, Lubich et al. \cite{Lubich2021} developed and analyzed high-order BDF schemes for the LLG equation. Unconditional unique solvability for such semi-implicit schemes has been established in \cite{jingrun2019analysis,Lubich2021}; however, their convergence analysis typically requires the temporal step size to be proportional to the spatial grid size. Despite these developments, the practical implementation and evaluation of high-order numerical methods for micromagnetic modeling continue to be active research areas. 
A notable limitation persists: 
most existing methods are limited to projection methods to handle the norm preserving constraints. This projection will preserve the magnetization direction, but with normalized magnitude. Such a projection step due to its nonlinearity bring a huge challenge in stability and convergence analysis.
% Furthermore, there is a scarcity of third-order accurate methods tailored for practical micromagnetic simulations, especially those capable of handling arbitrary damping parameters and facilitating quantitative comparisons among different numerical approaches. 
To address this gap, this paper presents a naive first-order accurate numerical scheme for solving the LLG equation with an existing method as a stable solution and a Crank-Nicolson's type predictor to preserve the normalizing constraint. Therefore, there's no constraints applied, just using the model itself to get a stable results. We also conduct comprehensive computational tests to verify the numerical stability of the proposed method. These tests further reveal that the accuracy and stability. 
% the dissipation properties of the new scheme differ from those of our earlier method \cite{xie2025schemeB} when applied to systems with large damping parameters.

The rest of this paper is organized as follows. \cref{sec: numerical scheme} begins with a review of the micromagnetic model, followed by a detailed description of the proposed numerical scheme and a comparative discussion with first-order Gauss-Seidel projection method \cite{wang2001gauss}. \cref{sec:experiments} presents extensive numerical results, encompassing verification of temporal and spatial accuracy in one-dimensional (1D) and three-dimensional (3D) settings.
% assessments of computational efficiency (via comparisons with BDF1 and BDF2), stability analysis with respect to the damping parameter, and an investigation of domain wall velocity dependence on both damping and external magnetic field. 
Concluding remarks and potential future research directions are provided in \cref{sec:conclusions}.

%\section{Formalism\label{Sec.Form}}
\section{The physical model and governing equation}
\label{sec: numerical scheme}

%\subsection{Governing equation}

The Landau-Lifshitz-Gilbert (LLG) equation forms the fundamental basis of micromagnetics, providing a rigorous description of the spatiotemporal evolution of magnetization in ferromagnetic materials by incorporating two key physical phenomena: gyromagnetic precession and dissipative relaxation \cite{Landau1935On,Brown1963micromagnetics}. In nondimensional form, this governing equation is expressed as
\begin{align}\label{c1-large}
{\m}_t =-{\m}\times{\bm h}_{\text{eff}}-\alpha{\m}\times({\m}\times{\bm h}_{\text{eff}})
\end{align}
subject to the homogeneous Neumann boundary condition
\begin{equation}\label{boundary-large}
\frac{\partial{\m}}{\partial {\bm \nu}}\Big|_{\partial \Omega}=0,
\end{equation}
where \(\Omega \subset \mathbb{R}^d\) (\(d=1,2,3\)) represents the bounded domain of the ferromagnetic material, and \(\bm \nu\) is the unit outward normal vector on the boundary \(\partial \Omega\). This boundary condition ensures no magnetic surface charge, a physically appropriate assumption for isolated ferromagnetic systems.

The magnetization field \(\m: \Omega \to \mathbb{R}^3\) is a three-dimensional vector field satisfying the pointwise constraint \(|\m| = 1\),  stemming from the quantum mechanical alignment of electron spins in ferromagnets. The first term on the right-hand side of \cref{c1-large} describes gyromagnetic precession, where magnetic moments precess around the effective field \(\bm h_{\text{eff}}\). The second term represents dissipative relaxation, with \(\alpha > 0\) being the dimensionless Gilbert damping coefficient that governs the rate of energy dissipation into the lattice.

From the perspective of the Gibbs free energy functional, the effective field \(\bm h_{\text{eff}}\) is obtained as the functional derivative of the Gibbs free energy \(F[\m]\) with respect to the magnetization, i.e., \(\bm h_{\text{eff}} = -\delta F/\delta \m\). his functional incorporates all relevant energy contributions in ferromagnetic systems—exchange, anisotropy, magnetostatic (stray), and Zeeman energies—and is given by
\begin{equation}\label{LL-Energy}
F[\m] = \frac {\mu_0 M_s^2}{2} \left\{\int_\Omega \left( \epsilon|\nabla\m|^2 +
q\left(m_2^2 + m_3^2\right)
-2\h_e\cdot\m - \h_s\cdot\m \right)\mathrm{d}\x \right\} . 
\end{equation}
Here, \(\mu_0 = 4\pi \times 10^{-7}\, \text{H/m}\) is the vacuum permeability, \(M_s\) is the saturation magnetization, and \(\epsilon\) and \(q\) are dimensionless parameters defined subsequently. The vectors \(\h_e\) and \(\h_s\) denote the externally applied magnetic field and the stray field, respectively. For uniaxial ferromagnetic materials with a single easy axis, the effective field \(\bm h_{\text{eff}}\) decomposes into distinct physical components, yielding the explicit expression
\begin{align}
{\bm h}_{\text{eff}} =\epsilon\Delta\m-q(m_2\e_2+m_3\e_3)+\h_s+\h_e,
\end{align}
where \(\epsilon = C_{\text{ex}}/(\mu_0 M_s^2 L^2)\) and \(q = K_u/(\mu_0 M_s^2)\). Here, \(L\) is the characteristic length scale, \(C_{\text{ex}}\) is the exchange constant (controlling short-range spin alignment), and \(K_u\) is the uniaxial anisotropy constant (representing the energy cost for deviation from the easy axis). The unit vectors \(\e_2 = (0,1,0)\) and \(\e_3 = (0,0,1)\) efine the hard axes perpendicular to the uniaxial easy axis, and \(\Delta\) is the Laplacian operator in \(d\)-dimensional space.
For Permalloy (NiFe), a common soft ferromagnetic material in spintronics, standard parameter values from the literature are: \(C_{\text{ex}} = 1.3 \times 10^{-11}\, \text{J/m}\), \(K_u = 100\, \text{J/m}^3\), and \(M_s = 8.0 \times 10^5\, \text{A/m}\). The stray field \(\h_s\) originates from magnetic charge distributions at domain boundaries and material surfaces and is mathematically represented by the integral expression
\begin{align}\label{eqn:div}
{\h}_{\text{s}}=\frac{1}{4\pi}\nabla \int_{\Omega} \nabla\left( \frac{1}{|\x-\y|}\right)\cdot {\bm m}({\bm y})\,d{\bm y},
\end{align}
which is a formulation that exhibits long-range spatial correlations. A critical computational advancement for practical micromagnetic simulations is that for rectangular domains \(\Omega\), the evaluation of \(\h_s\) can be efficiently computed via the Fast Fourier Transform (FFT) \cite{Wang2000}, which reduces the asymptotic computational complexity from \(O(N^d)\) to \(O(N^d \log N)\) for \(d\)-dimensional grids, enabling large-scale simulations.

To facilitate numerical discretization, we introduce the composite source term
\begin{align}\label{eq-4}
\f=-Q(m_2\e_2+m_3\e_3)+\h_s+\h_e.
\end{align}
which aggregates the anisotropy, stray field, and external field contributions. Substituting this source term into \cref{c1-large}, the LLG equation is re-expressed as
\begin{align}\label{eq-5}
\m_t=-\m\times(\epsilon\Delta\m+\f)-\alpha\m\times(\m\times(\epsilon\Delta\m+\f).)
\end{align}
Leveraging the vector triple product identity \(\a \times ({\bm b} \times {\bm c}) = (\a \cdot {\bm c}){\bm b} - (\a \cdot {\bm b}){\bm c}\) and the pointwise constraint \(|\m| = 1\) (which implies \(\m \cdot \partial_t \m = 0\) via time differentiation), we simplify \cref{eq-5} to an equivalent formulation that is more amenable to stable numerical discretization:
\begin{equation}\label{eq-model}
\m_t=\alpha  (\epsilon\Delta\m+\f)+\alpha \left(\epsilon |\nabla \m|^2 -\m \cdot\f \right)\m-\m\times(\epsilon\Delta\m+\f).
\end{equation}

\section{Proposed method}

For the construction of the structure preserving method, we set $\f=0$, $\epsilon=1$ and $\alpha=0$ for \cref{eq-5}, we have the LLG equation below,
\begin{align}\label{eq-alpha-0}
\m_t=-\m\times\Delta\m.
\end{align}
If we consider the simple linear vectorial equation 
\begin{align}\label{eq-CN}
\m_t=-\m\times \a,
\end{align}
where $\a^T=(a_1,a_2,a_3)$ is a constant vector. We use the Crank-Nicolson method to \cref{eq-CN}, we have
\begin{align}\label{eq-CN_1}
	\frac{\m_h^{n+1}-\m_h^n}{\Delta t}=-\frac{\m_h^{n+1}+\m_h^n}{2} \times \a,
\end{align}
which is norm preserving, since that if $\m_h^{n+1}+\m_h^n$ do the inner product for both sides, leading to
\begin{align*}
	\|\m_h^{n+1}\|_2=\|\m_h^n\|_2.
\end{align*}
After rearranging the compact form for \cref{eq-CN_1}, we have
\begin{align*}
	\begin{pmatrix}
	1&\frac12 \Delta t a_3&-\frac12 \Delta t a_2\\
	-\frac12 \Delta ta_3&1&\frac12 \Delta ta_1\\
	\frac12 \Delta ta_2&-\frac12 \Delta t a_1&1
	\end{pmatrix}\begin{pmatrix}
	m_1^{n+1}\\
	m_2^{n+1}\\
	m_3^{n+1}
	\end{pmatrix}=\begin{pmatrix}
	m_1^n+\frac12 \Delta t(a_2 m_3^n-a_3 m_2^n)\\
	m_2^n+\frac12 \Delta t(a_3 m_1^n-a_1 m_3^n)\\
	m_3^n+\frac12 \Delta t(a_1m_2^n-a_2 m_1^n)
	\end{pmatrix}
\end{align*}
We then have another form
\begin{align*}
	\begin{pmatrix}
	m_1^{n+1}\\
	m_2^{n+1}\\
	m_3^{n+1}
	\end{pmatrix}=\begin{pmatrix}
	1&\frac12 \Delta t a_3&-\frac12 \Delta t a_2\\
	-\frac12 \Delta ta_3&1&\frac12 \Delta ta_1\\
	\frac12 \Delta ta_2&-\frac12 \Delta t a_1&1
	\end{pmatrix}^{-1}\begin{pmatrix}
	1&-\frac12 \Delta t a_3 & \frac12 \Delta t a_3\\
	\frac12 \Delta t a_3&1& -\frac12 \Delta t a_1\\
	-\frac12 \Delta t a_2&\frac12 \Delta t a_1&1
	\end{pmatrix}\begin{pmatrix}
	m_1^n\\
	m_2^n\\
	m_3^n
	\end{pmatrix}=A\begin{pmatrix}
	m_1^n\\
	m_2^n\\
	m_3^n
	\end{pmatrix}
\end{align*}
where 
\begin{align*}
	A=\frac{1}{S}\begin{pmatrix}
	1+\beta^{2}a_1^{2}&-2\beta a_3+\beta^{2}a_1a_2&2\beta a_2+\beta^{2}a_1a_3\\
	2\beta a_3+\beta^{2}a_1a_2&1+\beta^{2}a_2^{2}&-2\beta a_1+\beta^{2}a_2a_3\\
	-2\beta a_2+\beta^{2}a_1a_3&2\beta a_1+\beta^{2}a_2a_3&1+\beta^{2}a_3^{2}
	\end{pmatrix}
\end{align*}
where $S=\det(A)=1+\beta^2(a_1^2+a_2^2+a_3^2)$ and $\beta=\frac{\Delta t}{2}$.\\
From \cref{eq-CN_1}, we can get
\begin{align}
m_h^{n+1}-m_h^{n} &= -\beta\,(m_h^{n+1}+m_h^{n})\times \mathbf{a}.
\end{align}
Suppose C is a skew-symmetric matrix, then we know
\begin{align*}
(m_h^{n+1}+m_h^{n})\times \mathbf{a} &= C(m_h^{n+1}+m_h^{n}).
\end{align*}
So we can get 
\begin{align}
m_h^{n+1}-m_h^{n} &= -\beta\,C(m_h^{n+1}+m_h^{n}).
\end{align}
By rearranging the terms,we know
\begin{align}
(I+\beta C)\,m^{n+1} &= (I-\beta C)\,m^{n}.
\end{align}
So we can get the form
\begin{align}
m^{n+1} = A m^n, \\
where\;A = (I+\beta C)^{-1}(I-\beta C).\label{eq:A_def}
\end{align}
Taking transpose of \cref{eq:A_def} and using $(M^{-1})^\top=(M^\top)^{-1}$,
\[
A^\top = (I-\beta C)^\top\bigl((I+\beta C)^{-1}\bigr)^\top
       = (I-\beta C)^\top\bigl((I+\beta C)^\top\bigr)^{-1}.
\]
Since $C^\top=-C$, we have $(I-\beta C)^\top=I+\beta C$ and $(I+\beta C)^\top=I-\beta C$, hence
\[
A^\top = (I+\beta C)(I-\beta C)^{-1}.
\]
Consequently,
\[
A^\top A
= (I+\beta C)(I-\beta C)^{-1}(I+\beta C)^{-1}(I-\beta C)
= I.
\]
Thus $A$ is an orthogonal matrix.
For any $x\in\mathbb{R}^3$,
\[
\|Ax\|_2^2 = (Ax)^\top(Ax)=x^\top(A^\top A)x=x^\top x=\|x\|_2^2,
\]
so $\|Ax\|_2=\|x\|_2$ for all $x$. By definition of the induced $2$-norm,
\[
\|A\|_2 = \max_{x\neq 0}\frac{\|Ax\|_2}{\|x\|_2}=1.
\]
Hence $\|A\|_2$ is bounded by $1$ and is independent of the time step $\Delta t$.

For \cref{eq-alpha-0}, we propose the following structure preserving schemes:
\begin{itemize}
	\item Sheme I: explicit for $\Delta \m$, 
	\begin{align}\label{eq-CN_2}
	\frac{\m_h^{n+1}-\m_h^n}{\Delta t}=-\frac{\m_h^{n+1}+\m_h^n}{2} \times \Delta_h \m_h^n,
	\end{align}
	which gives a CFL-type condition for the stability.
		\item Sheme II: implicit for $\Delta \m$, 
	\begin{align}\label{eq-CN_3}
	\frac{\m_h^{n+1}-\m_h^n}{\Delta t}=-\frac{\m_h^{n+1}+\m_h^n}{2} \times \Delta_h \m_h^{n+1},
	\end{align}
	which poses a difficulty for the nonlinear systems with high complexity.
	\item Scheme III: Semi-implicit method,
	\begin{align}\label{eq-CN_4}
	\frac{\m_h^{n+1}-\m_h^n}{\Delta t}=-\frac{\m_h^{n+1}+\m_h^n}{2} \times \Delta_h \tilde{\m}_h^{n+1}.
	\end{align}
	% where $\g_h^{s}=(I-\Delta t \Delta_h)^{-1}\m_h^s,\; s=n,n+1$. 
    To be specific, we propose
	\begin{align*}
	\begin{pmatrix}
	1&\frac12 \Delta t \Delta_h\tilde{m}_3^{n+1}&-\frac12 \Delta t \Delta_h\tilde{m}_2^{n+1}\\
	-\frac12 \Delta t\Delta_h\tilde{m}_3^{n+1}&1&\frac12 \Delta t\Delta_h\tilde{m}_1^{n+1}\\
	\frac12 \Delta t\Delta_h\tilde{m}_2^{n+1}&-\frac12 \Delta t\Delta_h \tilde{m}_1^{n+1}&1
	\end{pmatrix}\begin{pmatrix}
	m_1^{n+1}\\
	m_2^{n+1}\\
	m_3^{n+1}
	\end{pmatrix}=\begin{pmatrix}
	m_1^n+\frac12 \Delta t(\Delta_h\tilde{m}_2^{n+1} m_3^n-\Delta_h\tilde{m}_3^{n+1} m_2^n)\\
	m_2^n+\frac12 \Delta t(\Delta_h\tilde{m}_3^{n+1} m_1^n-\Delta_h\tilde{m}_1^{n+1} m_3^n)\\
	m_3^n+\frac12 \Delta t(\Delta_h\tilde{m}_1^{n+1}m_2^n-\Delta_h\tilde{m}_2^{n+1} m_1^n)
	\end{pmatrix}
	\end{align*}
	The trick is to handling the $\tilde{\m}_h^{n+1}$. If $\tilde{\m}_h^{n+1}=(I-k\Delta_h)^{-1}\m^n_h$, 
	% \begin{align*}
	% \begin{pmatrix}
	% 1&\frac12 \Delta t \Delta_hg_3^n&-\frac12 \Delta t \Delta_hg_2^n\\
	% -\frac12 \Delta t\Delta_hg_3^n&1&\frac12 \Delta t\Delta_hg_1^n\\
	% \frac12 \Delta t\Delta_hg_2^n&-\frac12 \Delta t \Delta_hg_1^n&1
	% \end{pmatrix}\begin{pmatrix}
	% m_1^{n+1}\\
	% m_2^{n+1}\\
	% m_3^{n+1}
	% \end{pmatrix}=\begin{pmatrix}
	% m_1^n+\frac12 \Delta t(\Delta_hg_2^n m_3^n-\Delta_hg_3^n m_2^n)\\
	% m_2^n+\frac12 \Delta t(\Delta_hg_3^n m_1^n-\Delta_hg_1^n m_3^n)\\
	% m_3^n+\frac12 \Delta t(\Delta_hg_1^n m_2^n-\Delta_hg_2^n m_1^n)
	% \end{pmatrix}
	% \end{align*}
	such a scheme is proved to be mildly better than CFL condition. We should get a much better estimation for $\tilde{\m}_h^{n+1}$.
\end{itemize}

The Gauss-Seidel projection method by \cite{wang2001gauss} is given below,
\begin{align*}
    m_1^{n+1}&=m_1^n+\Delta t \Delta_h g_2^n \cdot m_3^n-\Delta t \Delta_h g_3^n\cdot m_2^n,\\
    m_2^{n+1}&=m_2^n +\Delta t \Delta_h g_3^n\cdot m_1^{n+1}-\Delta t \Delta_h g_1^{n+1}\cdot m_3^n\\
    m_3^{n+1}&=m_3^n +\Delta t \Delta_h g_1^{n+1}\cdot m_2^{n+1}-\Delta
     t \Delta_h g_2^{n+1}\cdot m_1^{n+1},
\end{align*}
then with a projection step $\m_h^{n+1}=\frac{\m_h^{n+1}}{|\m_h^{n+1}|}$. Note that here we have
\begin{align*}
    \frac{m_i^{n+1}-m_i^n}{\Delta t} =\Delta_h m_i^{n+1},
\end{align*}
which gives $m_i^{n+1}=(I-\Delta
 t \Delta_h)^{-1}m_i^n:=g_i^n$ and  
 \begin{align*}
     (\m_h^{n+1}-\m_h^n)\times \m_h^n=\m_h^{n+1}\times \m_h^n=\g_h^{n}\times \m_h^n=\Delta t \Delta_h \m_h^{n+1}\times \m_h^n=\Delta t \Delta_h \g_h^{n+1}\times \m_h^n.
 \end{align*}
 Then the GSPM becomes
 \begin{align*}
    m_1^{n+1}&=m_1^n+g_2^n \cdot m_3^n- g_3^n\cdot m_2^n,\\
    m_2^{n+1}&=m_2^n + g_3^n\cdot m_1^{n+1}- g_1^{n+1}\cdot m_3^n\\
    m_3^{n+1}&=m_3^n + g_1^{n+1}\cdot m_2^{n+1}- g_2^{n+1}\cdot m_1^{n+1},
\end{align*}

Based on the stability, we propose a structure preserving method A for LLG equation without damping term below,
\begin{itemize}
    \item Gauss-Seidel iteration:
    \begin{align*}
    \tilde{m}_1^{n+1}&=m_1^n+g_2^n \cdot m_3^n- g_3^n\cdot m_2^n,\\
    \tilde{m}_2^{n+1}&=m_2^n + g_3^n\cdot \tilde{m}_1^{n+1}- \tilde{m}_1^{n+1}\cdot m_3^n\\
    \tilde{m}_3^{n+1}&=m_3^n + \tilde{g}_1^{n+1}\cdot \tilde{m}_2^{n+1}- \tilde{g}_2^{n+1}\cdot \tilde{m}_1^{n+1},
\end{align*}
\item preserving iteration:
\begin{align*}
	\begin{pmatrix}
	1&\frac12 \Delta t \Delta_h\tilde{m}_3^{n+1}&-\frac12 \Delta t \Delta_h\tilde{m}_2^{n+1}\\
	-\frac12 \Delta t\Delta_h\tilde{m}_3^{n+1}&1&\frac12 \Delta t\Delta_h\tilde{m}_1^{n+1}\\
	\frac12 \Delta t\Delta_h\tilde{m}_2^{n+1}&-\frac12 \Delta t\Delta_h \tilde{m}_1^{n+1}&1
	\end{pmatrix}\begin{pmatrix}
	m_1^{n+1}\\
	m_2^{n+1}\\
	m_3^{n+1}
	\end{pmatrix}=\begin{pmatrix}
	m_1^n+\frac12 \Delta t(\Delta_h\tilde{m}_2^{n+1} m_3^n-\Delta_h\tilde{m}_3^{n+1} m_2^n)\\
	m_2^n+\frac12 \Delta t(\Delta_h\tilde{m}_3^{n+1} m_1^n-\Delta_h\tilde{m}_1^{n+1} m_3^n)\\
	m_3^n+\frac12 \Delta t(\Delta_h\tilde{m}_1^{n+1}m_2^n-\Delta_h\tilde{m}_2^{n+1} m_1^n)
	\end{pmatrix}
	\end{align*}
\end{itemize}
We choose an exact solution as shown in the \cref{sec:experiments} in 1D, we get the result which is presented in \cref{tab-1}. 
\begin{table}[htbp]
	\centering
	\caption{Scheme A when $h = 5D-4$, $T=1d-1$.}\label{tab-1}
	\begin{tabular}{|c|c|c|c|}
		\hline
		$k$ & $\|\m_h-\m_e\|_\infty$ & $\|\m_h-\m_e\|_2$ &$\|\m_h-\m_e\|_{H^1}$ \\
		\hline
		2.0D-2 & 1.912933858544919&0.221808000745059&2.464765038269722e+02\\
		1.0D-2 &2.000796074933401&1.945760812419760&5.107532601677993e+02\\
		5.0D-3 &1.994436895670912&1.778123563276698&6.113345097672144e+02\\
		2.5D-3 &2.034119601928516&1.276682587287282&4.811857712154192e+02\\
		1.25D-3 &2.014352127693102&0.997384512352806&5.708736422654824e+02\\
		6.25D-4 &2.054232763749968&1.316576304183331&2.764074721845014e+02 \\
		\hline
	\end{tabular}
\end{table}

In turn, we propose a new structure preserving method B for LLG equation without damping term below,
\begin{itemize}
    \item Gauss-Seidel iteration:
    \begin{align*}
    \tilde{m}_1^{n+1}&=m_1^n+g_2^n \cdot m_3^n- g_3^n\cdot m_2^n,\\
    \tilde{m}_2^{n+1}&=m_2^n + g_3^n\cdot \tilde{m}_1^{n+1}- \tilde{m}_1^{n+1}\cdot m_3^n\\
    \tilde{m}_3^{n+1}&=m_3^n + \tilde{g}_1^{n+1}\cdot \tilde{m}_2^{n+1}- \tilde{g}_2^{n+1}\cdot \tilde{m}_1^{n+1},
\end{align*}
\item double diffusion iteration:
	\begin{align*}
		m_i^{*,n+1}=(I-\Delta t \Delta_h)^{-1} \tilde{{m}}_i^{n+1}.
	\end{align*}
\item preserving iteration:
\begin{align*}
	\begin{pmatrix}
	1&\frac12 \Delta t \Delta_h {m}_3^{*,n+1}&-\frac12 \Delta t \Delta_h{m}_2^{*,n+1}\\
	-\frac12 \Delta t\Delta_h{m}_3^{*,n+1}&1&\frac12 \Delta t\Delta_h{m}_1^{*,n+1}\\
	\frac12 \Delta t\Delta_h{m}_2^{*,n+1}&-\frac12 \Delta t\Delta_h {m}_1^{*,n+1}&1
	\end{pmatrix}\begin{pmatrix}
	m_1^{n+1}\\
	m_2^{n+1}\\
	m_3^{n+1}
	\end{pmatrix}=\begin{pmatrix}
	m_1^n+\frac12 \Delta t(\Delta_h\tilde{m}_2^{n+1} m_3^n-\Delta_h\tilde{m}_3^{n+1} m_2^n)\\
	m_2^n+\frac12 \Delta t(\Delta_h\tilde{m}_3^{n+1} m_1^n-\Delta_h\tilde{m}_1^{n+1} m_3^n)\\
	m_3^n+\frac12 \Delta t(\Delta_h\tilde{m}_1^{n+1}m_2^n-\Delta_h\tilde{m}_2^{n+1} m_1^n)
	\end{pmatrix}
	\end{align*}
\end{itemize}
We choose an exact solution as shown in the \cref{sec:experiments} in 1D, we get the result which is presented in \cref{tab-2} for the temporal accuracy and the spatial accuracy in \cref{tab-3}. 
\begin{table}[htbp]
	\centering
	\caption{Scheme B without damping when $h = 5D-4$, $T=1d-1$ for the temporal accuracy test in 1D.}\label{tab-2}
	\begin{tabular}{|c|c|c|c|}
		\hline
		$k$ & $\|\m_h-\m_e\|_\infty$ & $\|\m_h-\m_e\|_2$ &$\|\m_h-\m_e\|_{H^1}$ \\
		\hline
			2.0D-2 & 0.013537707345784&0.008843496854506&0.060208506436330\\
		1.0D-2 &0.009123251664833&0.005420865663023&0.039213187891314\\
		5.0D-3 &0.005104857017801&0.003077981908226&0.022450427816651\\
		2.5D-3 &0.002619628411882&0.001650564604991&0.011951250416275\\
		1.25D-3 &0.001312495917216&8.571574701860828e-04&0.006160144072687\\
		6.25D-4 &6.544134217098541e-04&4.377049821189154e-04&0.003134030652176 \\
		3.125D-4 &3.261646570654059e-04&2.214333020304648e-04&0.001583901331159 \\
		\hline
	\end{tabular}
\end{table}

\begin{table}[htbp]
	\centering
	\caption{Scheme B without damping when $k= 1D-5$, $T=1d-1$ for the spatial accuracy test in 1D.}\label{tab-3}
	\begin{tabular}{|c|c|c|c|}
		\hline
		$h$ & $\|\m_h-\m_e\|_\infty$ & $\|\m_h-\m_e\|_2$ &$\|\m_h-\m_e\|_{H^1}$  \\
		\hline
		1/16 &4.278398818205395e-04&2.963275220020730e-04&0.002202811741580\\
		1/24 &1.919211673657578e-04&1.341642420580365e-04&9.693232742511687e-04\\
		1/32 &1.090341965380506e-04&7.765516858659106e-05&5.412047553905197e-04\\
		1/48 &4.991127543831075e-05&3.749181553669659e-05&2.384007332139399e-04\\
		1/64 &3.127374664693705e-05&2.360257701529603e-05&1.355561887947331e-04\\
		\hline
	\end{tabular}
\end{table}

We choose an initial condition and set the source term to be zeros as shown in the \cref{sec:experiments} in 1D, we get the result for the norm preserving test which is presented in \cref{tab-4}

\begin{table}[htbp]
	\centering
	\caption{Scheme B without damping when $h = 5D-4$, $T=1d-1$ for the norm preserving test in 1D.}\label{tab-4}
	\begin{tabular}{|c|c|}
		\hline
		$k$ & $\|\|\m_h\|_2-1\|_\infty$ \\
		\hline
			2.0D-2 &1.332267629550188e-15 \\
		1.0D-2 &1.665334536937735e-15\\
		5.0D-3 &3.219646771412954e-15\\
		2.5D-3 &3.552713678800501e-15\\
		1.25D-3 &5.218048215738236e-15\\
		6.25D-4 &7.216449660063518e-15 \\
		3.125D-4 &1.354472090042691e-14 \\
		\hline
	\end{tabular}
\end{table}

For the model with damping term below,
\begin{align}\label{eq-alpha}
\m_t=-\m\times\Delta\m-\alpha \m \times (\m\times\Delta\m).
\end{align}
We consider the construction below,
	\begin{align}\label{eq-CN_4-1}
\frac{\m_h^{n+1}-\m_h^n}{\Delta t}=-\frac{\m_h^{n+1}+\m_h^n}{2} \times \Delta_h \tilde{\m}_h^{n+1}-\alpha\frac{\m_h^{n+1}+\m_h^n}{2} \times \left(\m_h^n\times \Delta_h \tilde{\m}_h^{n+1}\right),
\end{align}
which is also a structure preserving method with $\|\m_h^{n+1}\|_2=\|\m_h^{n}\|_2$. For such case, we propose the structure preserving method A below,
\begin{itemize}
    \item Gauss-Seidel iteration:
    \begin{align*}
    \tilde{m}_1^{n+1}&=m_1^n+g_2^n \cdot m_3^n- g_3^n\cdot m_2^n,\\
    \tilde{m}_2^{n+1}&=m_2^n + g_3^n\cdot \tilde{m}_1^{n+1}- \tilde{m}_1^{n+1}\cdot m_3^n\\
    \tilde{m}_3^{n+1}&=m_3^n + \tilde{g}_1^{n+1}\cdot \tilde{m}_2^{n+1}- \tilde{g}_2^{n+1}\cdot \tilde{m}_1^{n+1},
\end{align*}
\item damping iteration:
\begin{align*}
    m_i^{*,n+1}=\tilde{m}_i^{n+1}+a\alpha \Delta t\Delta_h {m}_i^{*,n+1}
\end{align*}
where $i=1,2,3$.
\item preserving iteration:
\begin{align*}
	% \begin{pmatrix}
	% 1&\frac12 \Delta t \Delta_h\tilde{m}_3^{n+1}&-\frac12 \Delta t \Delta_h\tilde{m}_2^{n+1}\\
	% -\frac12 \Delta t\Delta_h\tilde{m}_3^{n+1}&1&\frac12 \Delta t\Delta_h\tilde{m}_1^{n+1}\\
	% \frac12 \Delta t\Delta_h\tilde{m}_2^{n+1}&-\frac12 \Delta t\Delta_h \tilde{m}_1^{n+1}&1
	% \end{pmatrix}\begin{pmatrix}
	% m_1^{n+1}\\
	% m_2^{n+1}\\
	% m_3^{n+1}
	% \end{pmatrix}=\begin{pmatrix}
	% m_1^n+\frac12 \Delta t(\Delta_h\tilde{m}_2^{n+1} m_3^n-\Delta_h\tilde{m}_3^{n+1} m_2^n)\\
	% m_2^n+\frac12 \Delta t(\Delta_h\tilde{m}_3^{n+1} m_1^n-\Delta_h\tilde{m}_1^{n+1} m_3^n)\\
	% m_3^n+\frac12 \Delta t(\Delta_h\tilde{m}_1^{n+1}m_2^n-\Delta_h\tilde{m}_2^{n+1} m_1^n)
	% \end{pmatrix}
    \frac{\m_h^{n+1}-\m_h^n}{\Delta t}=-\frac{\m_h^{n+1}+\m_h^n}{2} \times \Delta_h{\m}_h^{*,n+1}-\alpha\frac{\m_h^{n+1}+\m_h^n}{2} \times \left(\m_h^n\times \Delta_h {\m}_h^{*,n+1}\right),
	\end{align*}
\end{itemize}

We choose an exact solution for the damping case as shown in the \cref{sec:experiments} in 1D, we get the result which is presented in \cref{tab-5}. 
\begin{table}[htbp]
	\centering
	\caption{Scheme A with damping $\alpha=0.01$ when $h = 5D-4$, $T=1d-1$.}\label{tab-5}
	\begin{tabular}{|c|c|c|c|}
		\hline
		$k$ & $\|\m_h-\m_e\|_\infty$ & $\|\m_h-\m_e\|_2$ &$\|\m_h-\m_e\|_{H^1}$ \\
		\hline
		2.0D-2 & 0.032270043090955& 0.019087832991816&0.279301695885164\\
		1.0D-2 & 1.803719315427388& 1.600224788803020& 1.340508110980241e+03\\
		5.0D-3 & 1.998107768985083& 1.414039591515378& 1.422501264852177e+03\\
		2.5D-3 & 1.988093537073624& 1.413849232401981& 1.437369582244364e+03\\
		1.25D-3 & 2.004058157997310& 1.414275131753937& 1.342262038804035e+03\\
		6.25D-4 & 2.019779868128165& 1.414291912266987& 1.013529493054862e+03 \\
		\hline
	\end{tabular}
\end{table}

\begin{table}[htbp]
	\centering
	\caption{The temporal accuracy and spatial accuracy for the proposed method A with damping $\alpha=0.01$, $T=0.1$ with $k=h^2$ in 3D. }\label{tab-a-10-time-3D-Q-2}
	\begin{tabular}{|c|c|c|c|c|}
		\hline
		$k$&$h$ & $\|\m_h-\m_e\|_\infty$ & $\|\m_h-\m_e\|_2$ &$\|\m_h-\m_e\|_{H^1}$ \\
	\hline
	T/10&1/10&4.999577379423137e-04&2.886005530405278e-04&3.342953283317248e-04 \\
	T/40 &1/20&2.023563574676621&0.842764687252767&8.745291288589728 \\
	T/57 &1/24&-&-&-\\
	% T/78&1/28&&&\\
	% T/102 &1/32&&&\\
	% T/129 &1/36& &&\\
	\hline 
	\end{tabular}
\end{table}

In turn, we propose a new structure preserving method B for LLG equation with damping term below,
\begin{itemize}
    \item Gauss-Seidel iteration:
    \begin{align*}
    \tilde{m}_1^{n+1}&=m_1^n+g_2^n \cdot m_3^n- g_3^n\cdot m_2^n,\\
    \tilde{m}_2^{n+1}&=m_2^n + g_3^n\cdot \tilde{m}_1^{n+1}- \tilde{m}_1^{n+1}\cdot m_3^n\\
    \tilde{m}_3^{n+1}&=m_3^n + \tilde{g}_1^{n+1}\cdot \tilde{m}_2^{n+1}- \tilde{g}_2^{n+1}\cdot \tilde{m}_1^{n+1},
\end{align*}
\item damping iteration:
\begin{align*}
    m_i^{*,n+1}=\tilde{m}_i^{n+1}+\alpha \Delta t\Delta_h {m}_i^{*,n+1}
\end{align*}
where $i=1,2,3$.
\item double diffusion iteration:
\begin{align*}
		m_i^{**,n+1}=(I-\Delta t \Delta_h)^{-1} {{m}}_i^{*,n+1}.
	\end{align*}
\item preserving iteration:
\begin{align*}
	% \begin{pmatrix}
	% 1&\frac12 \Delta t \Delta_h\tilde{m}_3^{n+1}&-\frac12 \Delta t \Delta_h\tilde{m}_2^{n+1}\\
	% -\frac12 \Delta t\Delta_h\tilde{m}_3^{n+1}&1&\frac12 \Delta t\Delta_h\tilde{m}_1^{n+1}\\
	% \frac12 \Delta t\Delta_h\tilde{m}_2^{n+1}&-\frac12 \Delta t\Delta_h \tilde{m}_1^{n+1}&1
	% \end{pmatrix}\begin{pmatrix}
	% m_1^{n+1}\\
	% m_2^{n+1}\\
	% m_3^{n+1}
	% \end{pmatrix}=\begin{pmatrix}
	% m_1^n+\frac12 \Delta t(\Delta_h\tilde{m}_2^{n+1} m_3^n-\Delta_h\tilde{m}_3^{n+1} m_2^n)\\
	% m_2^n+\frac12 \Delta t(\Delta_h\tilde{m}_3^{n+1} m_1^n-\Delta_h\tilde{m}_1^{n+1} m_3^n)\\
	% m_3^n+\frac12 \Delta t(\Delta_h\tilde{m}_1^{n+1}m_2^n-\Delta_h\tilde{m}_2^{n+1} m_1^n)
	% \end{pmatrix}
    \frac{\m_h^{n+1}-\m_h^n}{\Delta t}=-\frac{\m_h^{n+1}+\m_h^n}{2} \times \Delta_h {\m}_h^{**,n+1}-\alpha\frac{\m_h^{n+1}+\m_h^n}{2} \times \left(\m_h^n\times \Delta_h {\m}_h^{**,n+1}\right),
	\end{align*}
\end{itemize}
We choose an exact solution for the damping case as shown in the \cref{sec:experiments} in 1D, we get the result which is presented in \cref{tab-6} for the temporal accuracy and the spatial accuracy in \cref{tab-7}.

We choose an initial condition and set the source term to be zeros as shown in the \cref{sec:experiments} in 1D, we get the result for the norm preserving test which is presented in \cref{tab-8}

For the full LLG equation, we propose the structure preserving method as below,
\begin{itemize}
    \item Gauss-Seidel iteration:
    \begin{align*}
    \tilde{m}_1^{n+1}&=m_1^n+g_2^n \cdot m_3^n- g_3^n\cdot m_2^n,\\
    \tilde{m}_2^{n+1}&=m_2^n + g_3^n\cdot \tilde{m}_1^{n+1}- \tilde{m}_1^{n+1}\cdot m_3^n\\
    \tilde{m}_3^{n+1}&=m_3^n + \tilde{g}_1^{n+1}\cdot \tilde{m}_2^{n+1}- \tilde{g}_2^{n+1}\cdot \tilde{m}_1^{n+1},
\end{align*}
where $g_i^n=(I-\epsilon\Delta t \Delta_h)^{-1}(m_i^n+\Delta t f_i^n)$, and $\tilde{g}_i^n=(I-\epsilon\Delta t \Delta_h)^{-1}(\tilde{m}_i^n+\Delta t \tilde{f}_i^n)$
\item damping iteration:
\begin{align*}
    m_i^{*,n+1}=\tilde{m}_i^{n+1}+\alpha \Delta t(\epsilon \Delta_h {m}_i^{*,n+1}+\tilde{f}_i^n)
\end{align*}
where $i=1,2,3$.
\item double diffusion iteration:
\begin{align*}
		m_i^{**,n+1}=(I-\epsilon\Delta t \Delta_h)^{-1} ({{m}}_i^{*,n+1}+f_i^{*,n+1}).
	\end{align*}
\item preserving iteration:
\begin{align*}
    \frac{\m_h^{n+1}-\m_h^n}{\Delta t}=-\frac{\m_h^{n+1}+\m_h^n}{2} \times (\epsilon \Delta_h {\m}_h^{**,n+1}+\f_h^{**,n+1})-\alpha\frac{\m_h^{n+1}+\m_h^n}{2} \times \left(\m_h^n\times (\epsilon \Delta_h {\m}_h^{**,n+1}+\f_h^{**,n+1})\right).
	\end{align*}
\end{itemize}

\section{Numerical experiments}
\label{sec:experiments}

\subsection{Accuracy and efficiency tests}

To simplify the accuracy verification, we set parameter  \(\epsilon=1\) and forcing term \(\f=0\) in the governing model \cref{eq-5}. Analytical exact solutions are derived for both one-dimensional (1D) and three-dimensional (3D) scenarios to serve as benchmarks for error quantification.

For the 1D case, the exact magnetization solution \(\m_e\) is:
\begin{equation*}
\m_e=\left(\cos(\cos(\pi x))\sin t, \sin(\cos(\pi x))\sin t, \cos t\right)^T,
\end{equation*}
while the corresponding 3D exact solution is:
\begin{equation*}
\m_e=\left(\cos(XYZ)\sin t, \sin(XYZ)\sin t, \cos t\right)^T,
\end{equation*}
where $X=x^2(1-x)^2$, $Y=y^2(1-y)^2$ and $Z=z^2(1-z)^2$.

These exact solutions satisfy the governing equation \cref{eq-5} when the forcing term is defined as \(\f_e=\partial_t \m_e+\m_e \times \Delta \m_e+\alpha\m_e\times (\m_e \times \Delta \m_e)\). They also comply with the homogeneous Neumann boundary condition, ensuring consistency with simulation constraints.

To isolate the temporal approximation error from spatial discretization effects, the spatial resolution in the 1D test is fixed at \(h=5\times 10^{-4}\)—a sufficiently fine grid that makes spatial error negligible compared to temporal error. The Gilbert damping parameter is set to \(\alpha=0.01\), and simulations run until final time \(T=0.1\). Under this configuration, the measured error primarily reflects temporal discretization inaccuracy.

The 3D temporal accuracy test faces inherent constraints from spatial resolution, as excessively fine grids incur prohibitive computational cost. To balance spatial and temporal error contributions, we adopt a coordinated refinement strategy for spatial mesh sizes (\(h_x, h_y, h_z\)) and temporal step-size (\(\Delta t\)) tailored to the proposed method's order: \(\Delta t=h_x^2=h_y^2=h_z^2=h^2=T/N_0\). Here, \(N_0\) is a refinement level parameter, with specific values given in subsequent results. Consistent with the 1D test, \(\alpha=0.01\), and the final time \(T\) is specified later. The first order temporal accuracy is verified from \cref{tab-6} and \cref{tab-a-10-time-3D-Q-2} for 1D and 3D tests, respectively. 

\begin{table}[htbp]
	\centering
	\caption{Scheme B with damping $\alpha=0.01$ when $h = 5D-4$, $T=1d-1$ for the temporal accuracy test in 1D.}\label{tab-6}
	\begin{tabular}{|c|c|c|c|}
		\hline
		$k$ & $\|\m_h-\m_e\|_\infty$ & $\|\m_h-\m_e\|_2$ &$\|\m_h-\m_e\|_{H^1}$ \\
		\hline
	    2.0D-2 &0.013517320602480 &0.008817014136207&0.060170606240083\\
		1.0D-2 &0.009093911836419&0.005399671144500&0.039122839681984\\
		5.0D-3 &0.005080840309998&0.003060550005000&0.022350170425717\\
		2.5D-3 &0.002605046974160&0.001639087410385&0.011882958530157\\
		1.25D-3 &0.001304955694311&8.505492063928758e-04&0.006121659921789\\
		6.25D-4 &6.507998847741087e-04&4.341324733964456e-04&0.003113488007849 \\
		3.125D-4 &3.244944640112823e-04&2.195635788923884e-04&0.001573104885449 \\
        \hline
        order &0.918275397527218 &0.896549694749004 &0.890829368604443 \\
		\hline
	\end{tabular}
\end{table}
% \begin{table}[htbp]
% 	\centering
% 	\caption{The temporal accuracy for proposed method with damping $\alpha=0.01$when $h = 5D-4$, $T=1d-1$.}\label{tab-a-v-3}
% 	\begin{tabular}{|c|c|c|c|}
% 		\hline
% 		$k$ & $\|\m_h-\m_e\|_\infty$ & $\|\m_h-\m_e\|_2$ &$\|\m_h-\m_e\|_{H^1}$\\
% 		\hline
% 		2.0D-2 &&& \\
% 		1.0D-2 &&&\\
% 		5.0D-3 &&& \\
% 		2.5D-3 &&& \\
% 		1.25D-3 &&&\\
% 		6.25D-4 &&& \\
% 		3.125D-4 &&& \\
% 		\hline 
% 		order &&&\\
% 		\hline
% 	\end{tabular}
% \end{table}

\begin{table}[htbp]
	\centering
	\caption{Scheme B with damping $\alpha=0.01$ when $k= 1D-5$, $T=1d-1$ for the spatial accuracy test in 1D.}\label{tab-7}
	\begin{tabular}{|c|c|c|c|}
		\hline
		$h$ & $\|\m_h-\m_e\|_\infty$ & $\|\m_h-\m_e\|_2$ &$\|\m_h-\m_e\|_{H^1}$  \\
		\hline
		1/16 &4.244649934947095e-04&2.939188755381458e-04&0.002202247993444\\
		1/24 &1.904974803562040e-04&1.330580408930415e-04&9.690131243557937e-04\\
		1/32 &1.082529259365597e-04&7.701029722223597e-05&5.408447724250333e-04\\
		1/48 &4.947541044093839e-05&3.717722990660095e-05&2.379203616012192e-04\\
		1/64 &3.101638869329459e-05&2.340268587852941e-05&1.349769713496774e-04\\
		\hline
        order &1.900919204662317&1.830084876345314&2.016814249261991\\
        \hline
	\end{tabular}
\end{table}

% \begin{table}[htbp]
% 	\centering
% 	\caption{The spatial accuracy for proposed method with damping $\alpha=0.01$ when $k= 1D-6$, $T=1d-1$.}\label{tab-a-10-space-Q-1}
% 	\begin{tabular}{|c|c|c|c|}
% 		\hline
% 		$h$ & $\|\m_h-\m_e\|_\infty$ & $\|\m_h-\m_e\|_2$ &$\|\m_h-\m_e\|_{H^1}$ \\
% 		\hline
% 		1/16 &&&\\
% 		1/24 &&&\\
% 		1/32 &&&\\
% 		1/48 &&&\\
% 		1/64 &&&\\
% 		\hline 
% 		order &&&\\
% 		\hline
% 	\end{tabular}
% \end{table}

Following temporal accuracy evaluation, spatial accuracy tests were conducted to quantify the spatial discretization performance of the proposed method. To prevent temporal errors from interfering with spatial accuracy assessment, the temporal step size was fixed at a sufficiently small \(k=10^{-6}\) for 1D test shown in \cref{tab-7}—making temporal errors negligible compared to spatial errors. We take $k=h^2$, such that the second order spatial accuracy is observed, which is consistent with 1D test.

\begin{table}[htbp]
	\centering
	\caption{The temporal accuracy and spatial accuracy for the proposed method with damping $\alpha=0.01$, $T=0.1$ with $k=h^2$ in 3D. }\label{tab-a-10-time-3D-Q-2}
	\begin{tabular}{|c|c|c|c|c|}
		\hline
		$k$&$h$ & $\|\m_h-\m_e\|_\infty$ & $\|\m_h-\m_e\|_2$ &$\|\m_h-\m_e\|_{H^1}$ \\
	\hline
	T/10&1/10&5.000319903357697e-04&2.886017045040417e-04&3.310372227939050e-04 \iffalse 2.326619 (s) \fi \\
	T/40 &1/20&1.260826539929427e-04&7.235669771008058e-05&1.188710835928558e-04 \iffalse 36.394504 (s) \fi \\
	T/57 &1/24&8.888966024822587e-05&5.093893875864552e-05&9.792411451232800e-05 \iffalse 88.084476  (s) \fi\\
	T/78&1/28&6.531622705829854e-05&3.742972672128460e-05&8.478809887279647e-05 \iffalse 195.984342 (s) \fi \\
	T/102 &1/32&5.025462550745097e-05&2.886200671615646e-05&7.652231513783945e-05 \iffalse 432.664799  (s) \fi\\
	T/129 &1/36& 4.000667801529190e-05&2.309319543965180e-05&7.100714565453237e-05 \iffalse 908.984384 (s) \fi \\
	\hline 
	order &&0.988367654372780&0.989451989230767&0.612888163729896\\
		\hline
	&	order &1.972428431148607&1.974601917925088&1.223271316109047\\
		\hline
	\end{tabular}
\end{table}

\subsection{Norm preserving tests}

We choose the initial condition with below,
\begin{align*}
    \m_0=\left(\cos(\cos(\pi x))\sin (0.01), \sin(\cos(\pi x))\sin (0.01), \cos (0.01)\right)^T,
\end{align*}
and \begin{align*}
\m_0=\left(\cos(XYZ)\sin (0.01), \sin(XYZ)\sin (0.01), \cos (0.01)\right)^T,
\end{align*}
in 1D and 3D, respectively. The results are presented in \cref{tab-8} and \cref{tab-9}.

\begin{table}[htbp]
	\centering
	\caption{Scheme B with damping $\alpha=0.01$ when $h = 5D-4$, $T=1d-1$ for the norm preserving test in 1D.}\label{tab-8}
	\begin{tabular}{|c|c|}
		\hline
		$k$ & $\|\|\m_h\|_2-1\|_\infty$ \\
		\hline
	    2.0D-2 & 1.665334536937735e-15 \\
		1.0D-2 & 2.220446049250313e-15 \\
		5.0D-3 & 3.441691376337985e-15\\
		2.5D-3 & 3.996802888650564e-15 \\
		1.25D-3 & 5.884182030513330e-15 \\
		6.25D-4 & 8.215650382226158e-15 \\
		3.125D-4 &1.054711873393899e-14 \\
		\hline
	\end{tabular}
\end{table}

\begin{table}[htbp]
	\centering
	\caption{The norm preserving test for the proposed method with damping $\alpha=0.01$, $T=0.1$ with $k=h^2$ in 3D. }\label{tab-9}
	\begin{tabular}{|c|c|c|}
		\hline
		$k$&$h$ & $\|\|\m_h\|_2-1\|_\infty$  \\
	\hline
	T/10&1/10&8.963940700823514e-13\\
	T/40 &1/20&5.582201367815287e-13\\
	T/57 &1/24&5.104805467226470e-13 \\
	T/78&1/28& 4.438671652451376e-13\\
	% T/102 &1/32&\\
	% T/129 &1/36&  \\
	\hline 
	\end{tabular}
\end{table}

\subsection{Initial conditions and stability tests}

In this section, we choose different initial conditions specified in each test.

In 1D, we specify the initial condition as
\begin{align*}
    \m_0&=\left(\cos(\cos(\pi x))\sin (0), \sin(\cos(\pi x))\sin (0), \cos (0)\right)^T\\
    \m_0&=\left(\cos(\cos(\pi x))\sin (0.01), \sin(\cos(\pi x))\sin (0.01), \cos (0.01)\right)^T.
\end{align*}

\begin{figure}[htbp]
    \centering
    \subfloat[$m_1$]{\includegraphics[width=0.5\linewidth]{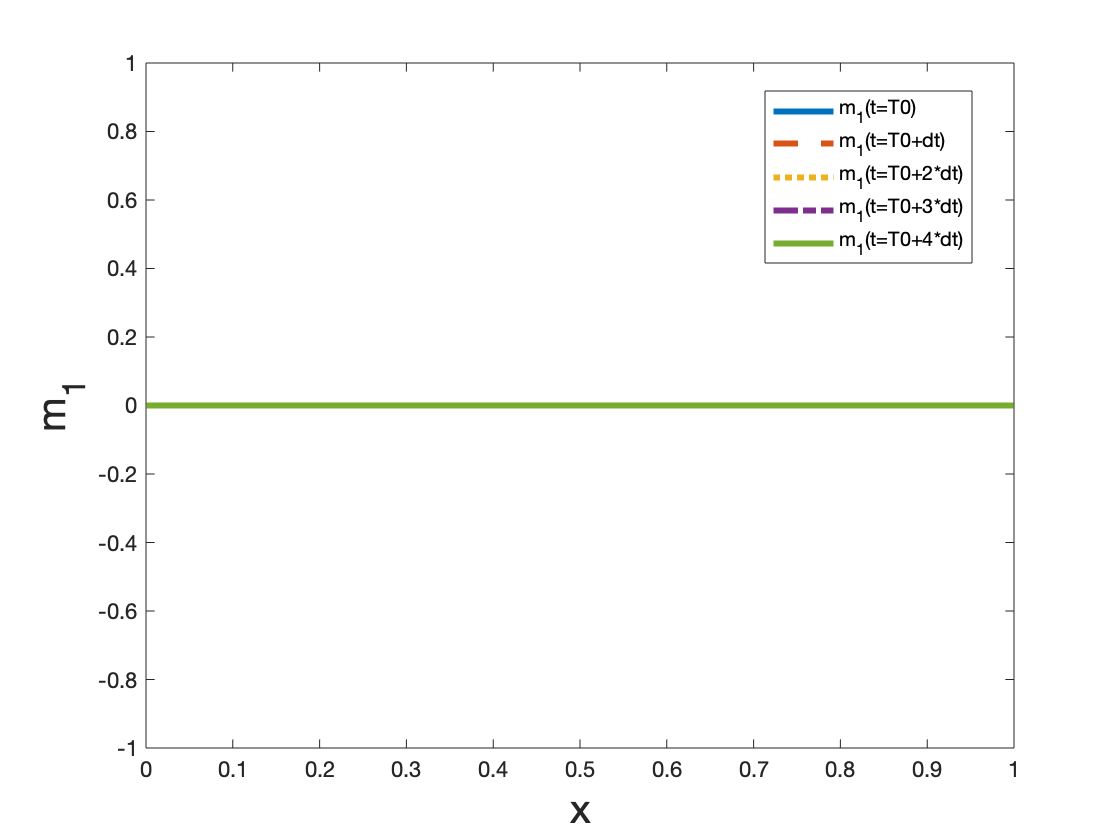}}
    \subfloat[$m_1$]{\includegraphics[width=0.5\linewidth]{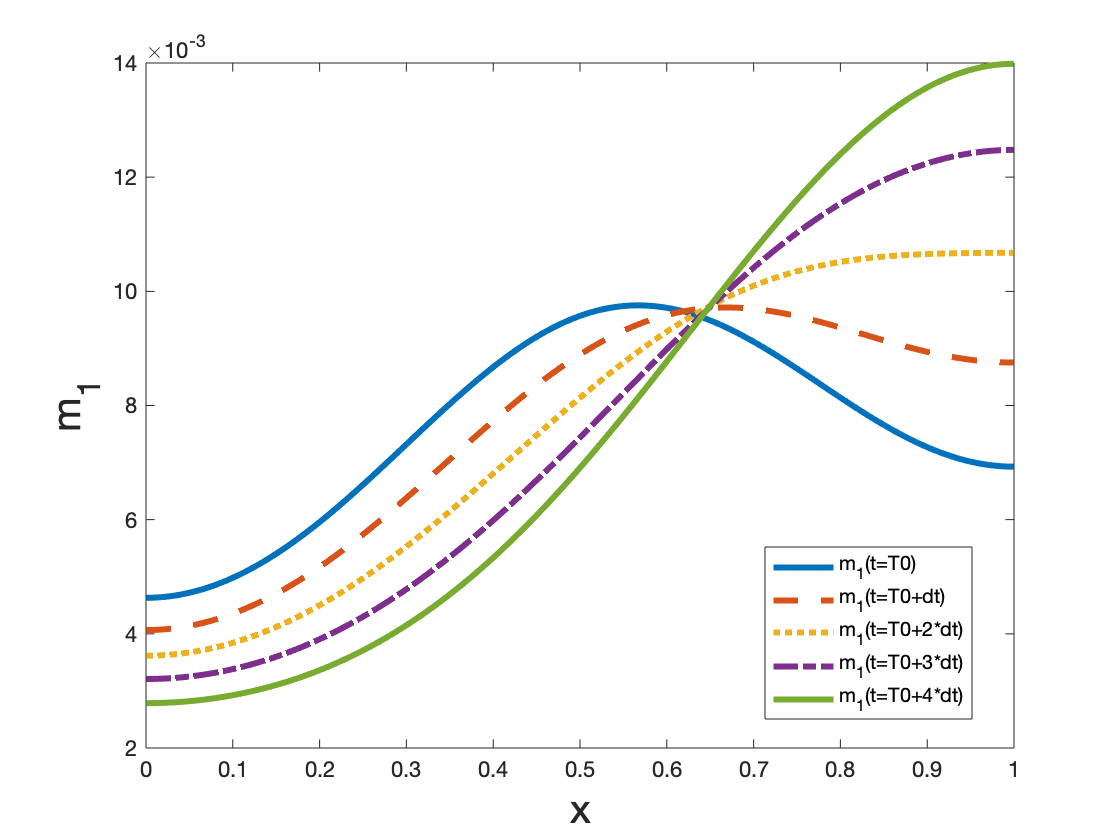}}
    \hspace{0.1in}
    \subfloat[$m_2$]{\includegraphics[width=0.5\linewidth]{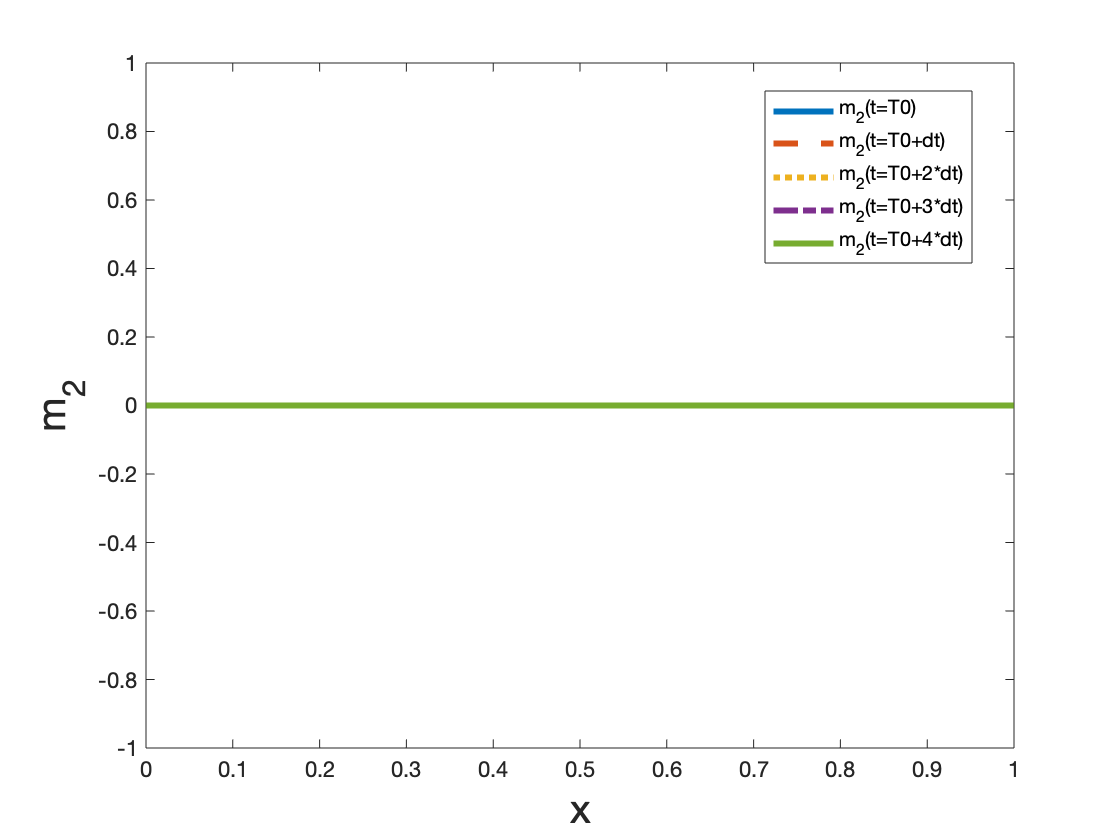}}
    \subfloat[$m_2$]{\includegraphics[width=0.5\linewidth]{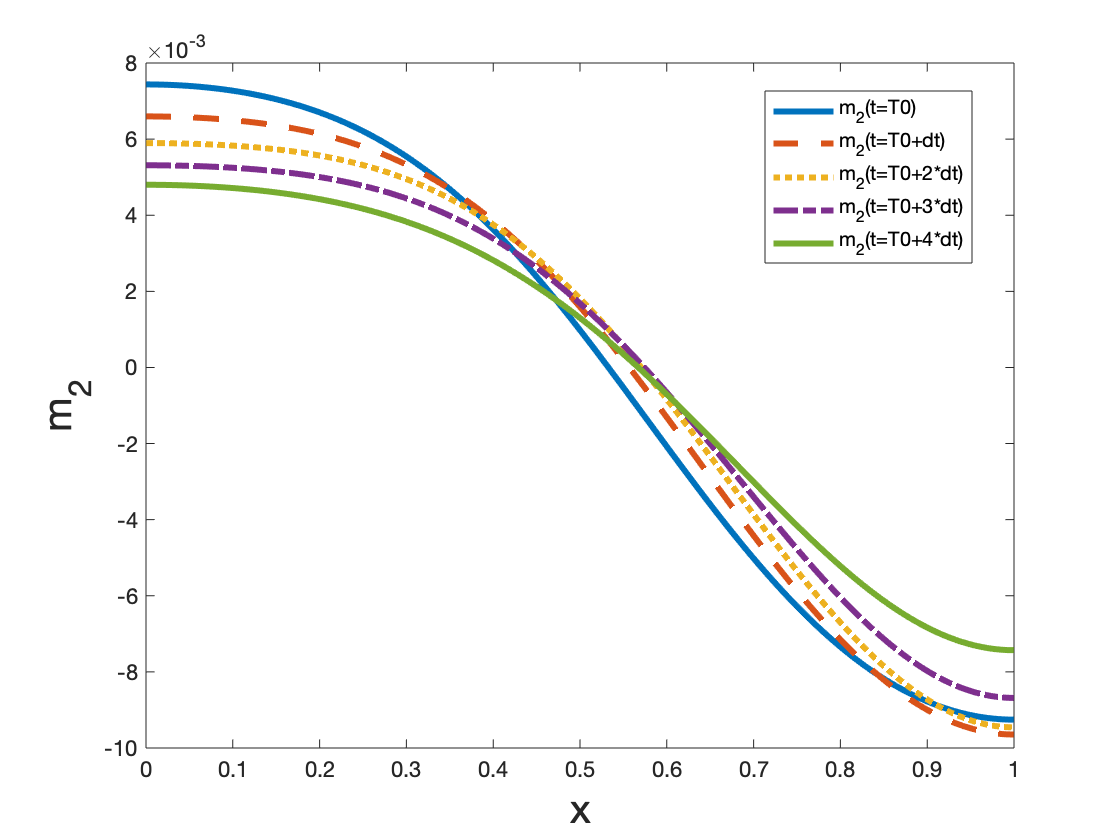}}
    \hspace{0.1in}
    \subfloat[$m_3$]{\includegraphics[width=0.5\linewidth]{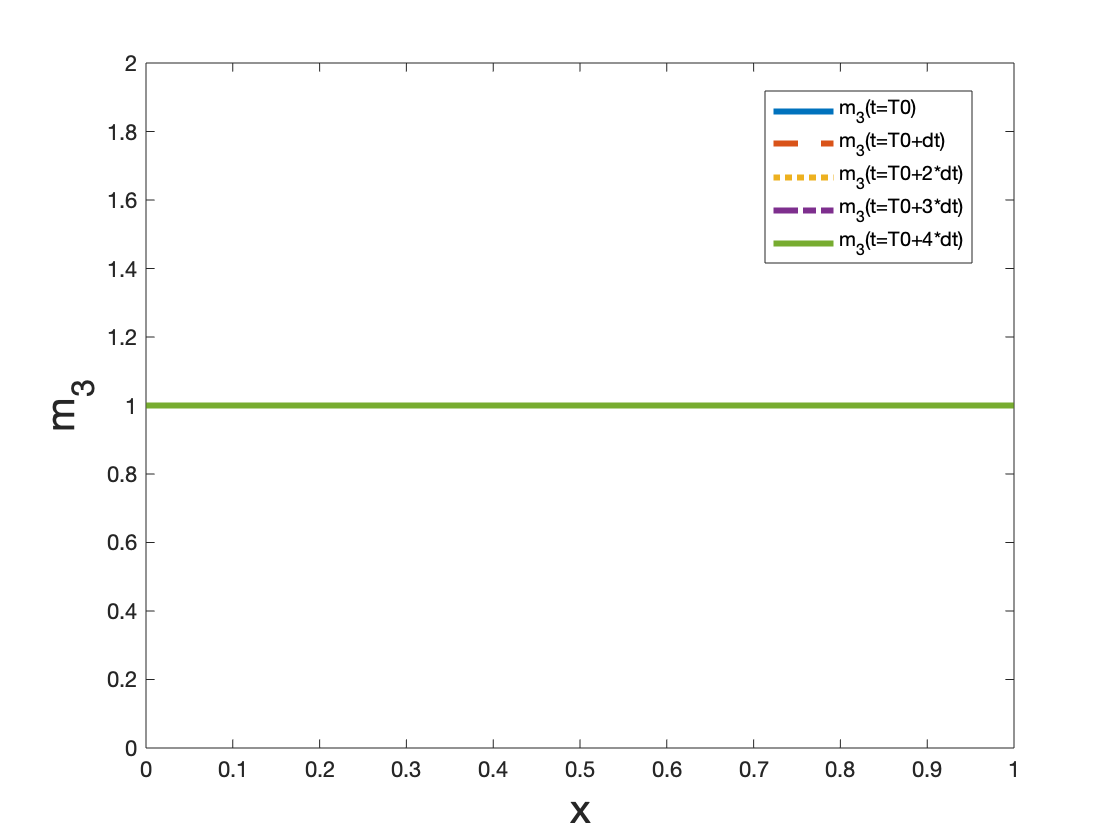}}
    \subfloat[$m_3$]{\includegraphics[width=0.5\linewidth]{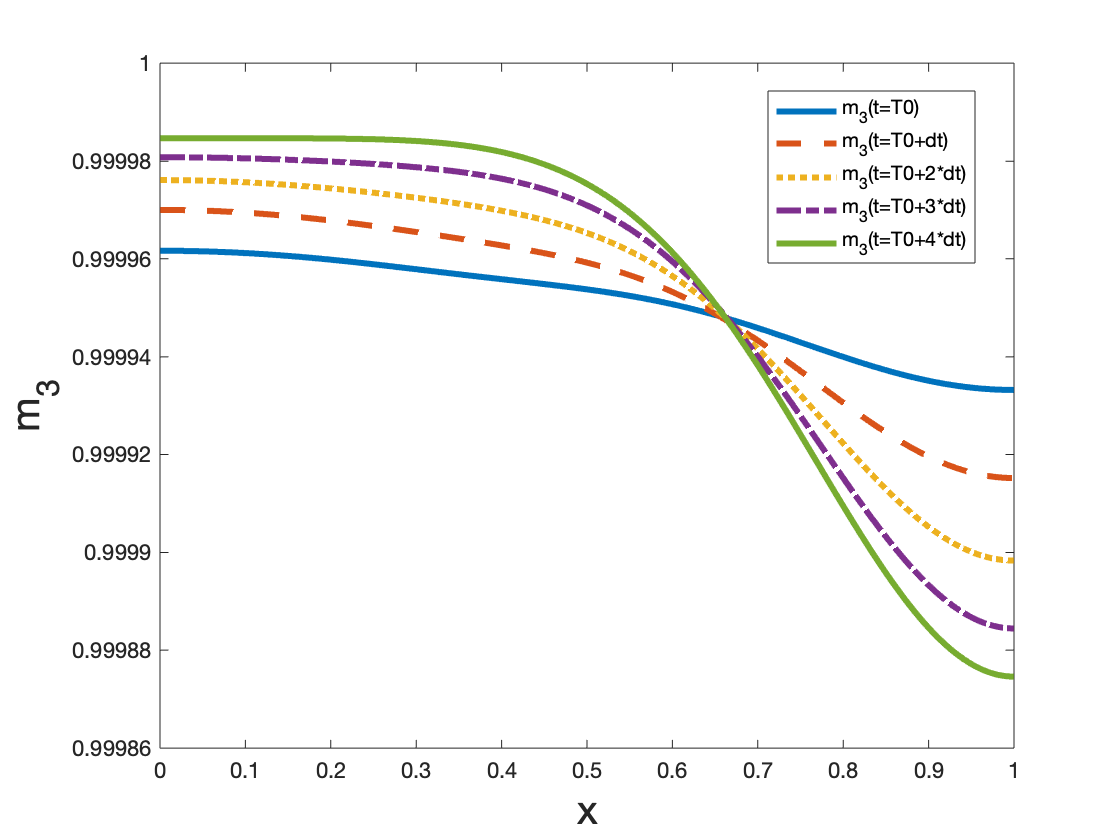}}
    \caption{The solution profile using proposed method in 1D given the initial condition $m_0$ without source term, $\alpha=0.01$ and $T=0.1$, $N_x=2000$, $N_t=5$. Left panel: the initial condition is given by $\m_0=\left(\cos(\cos(\pi x))\sin (0), \sin(\cos(\pi x))\sin (0), \cos (0)\right)^T$; Right panel: the initial condition is given by $\m_0=\left(\cos(\cos(\pi x))\sin (0.01), \sin(\cos(\pi x))\sin (0.01), \cos (0.01)\right)^T$.}
    \label{fig:1}
\end{figure}

\begin{figure}[htbp]
    \centering
    \subfloat[$m_1$]{\includegraphics[width=0.5\linewidth]{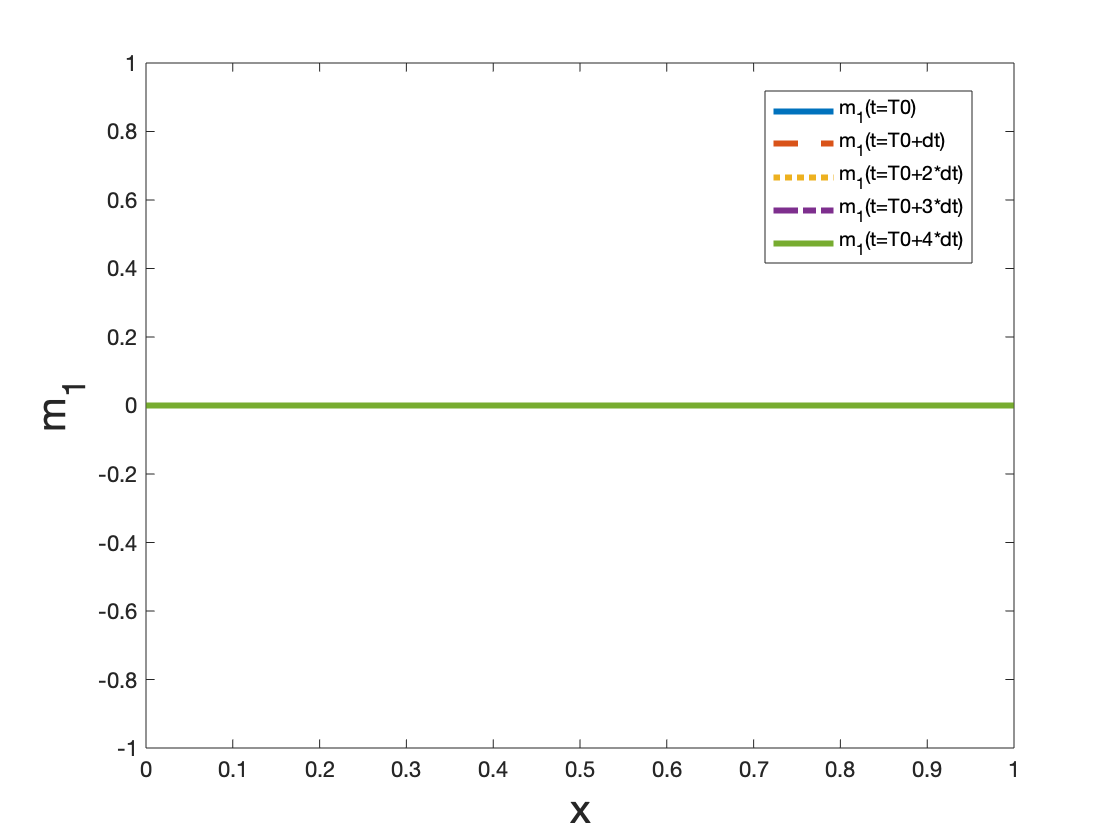}}
    \subfloat[$m_1$]{\includegraphics[width=0.5\linewidth]{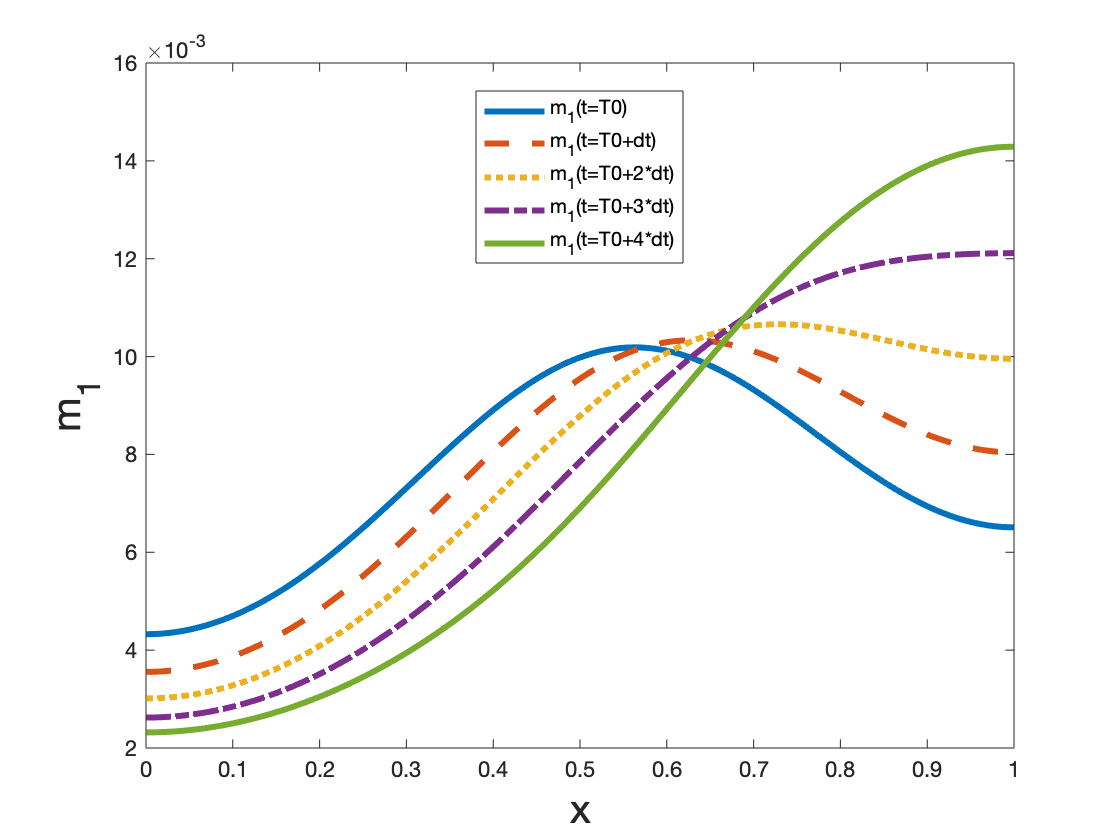}}
    \hspace{0.1in}
    \subfloat[$m_2$]{\includegraphics[width=0.5\linewidth]{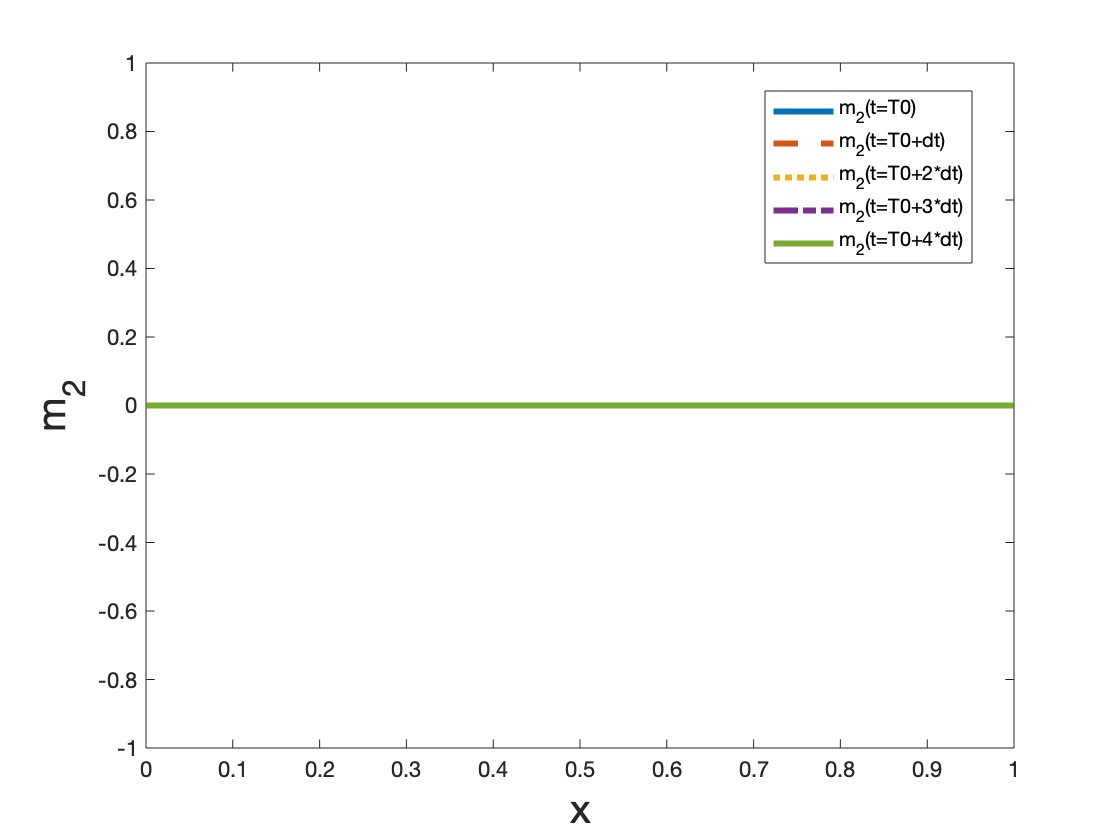}}
    \subfloat[$m_2$]{\includegraphics[width=0.5\linewidth]{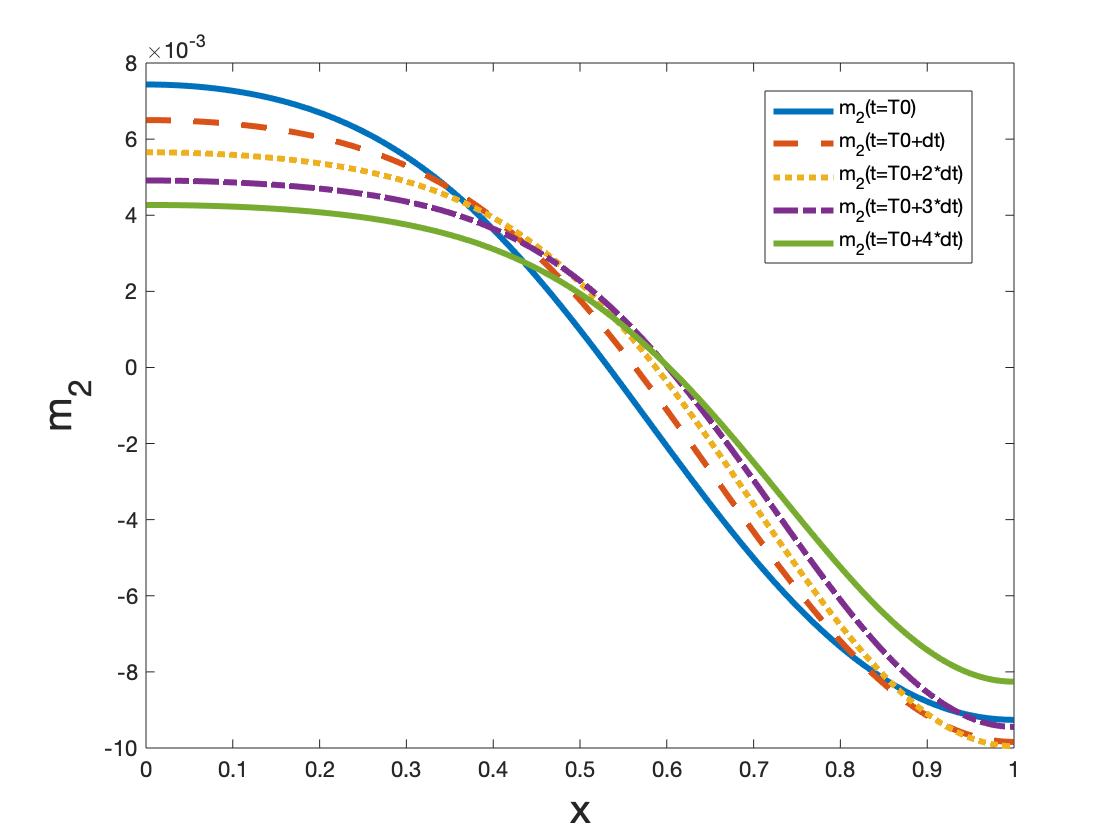}}
    \hspace{0.1in}
    \subfloat[$m_3$]{\includegraphics[width=0.5\linewidth]{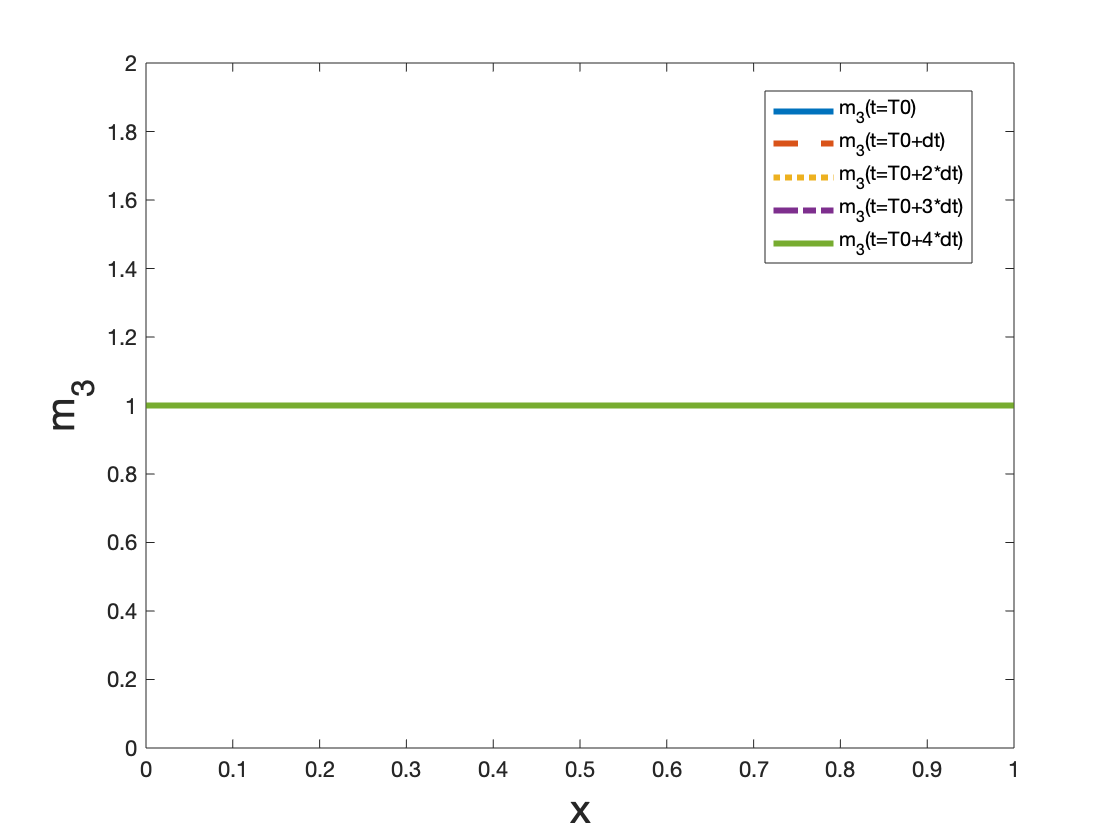}}
    \subfloat[$m_3$]{\includegraphics[width=0.5\linewidth]{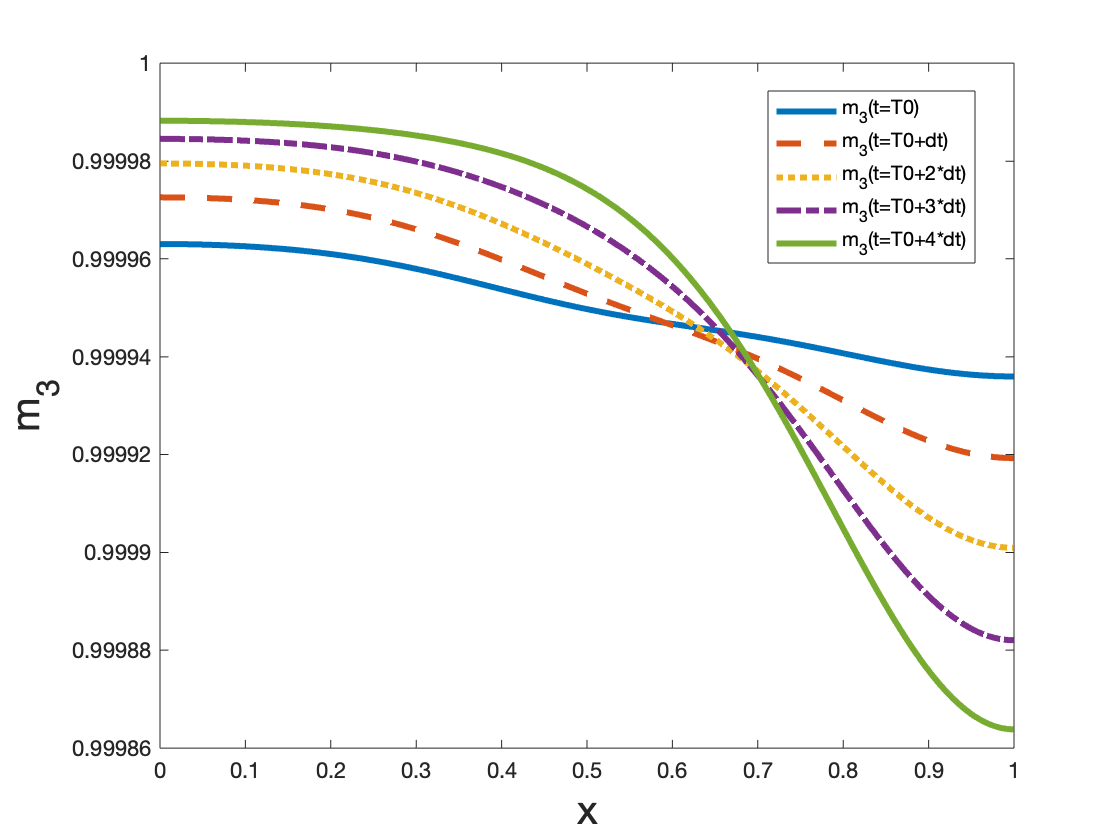}}
    \caption{The solution profile using GSPM method in 1D given the initial condition $m_0$ without source term, $\alpha=0.01$ and $T=0.1$, $N_x=2000$, $N_t=5$. Left panel: the initial condition is given by $\m_0=\left(\cos(\cos(\pi x))\sin (0), \sin(\cos(\pi x))\sin (0), \cos (0)\right)^T$; Right panel: the initial condition is given by $\m_0=\left(\cos(\cos(\pi x))\sin (0.01), \sin(\cos(\pi x))\sin (0.01), \cos (0.01)\right)^T$.}
    \label{fig:2}
\end{figure}

The results for 1D case using the proposed method and the first order GSPM to show the numerical consistency are presented in \Cref{fig:1} and \Cref{fig:2} with the paprameters $\alpha=0.01$, $N_x=2000$ and $N_t=5$. 

\begin{figure}[htbp]
    \centering
    \subfloat[arrow profile]{\includegraphics[width=0.5\linewidth]{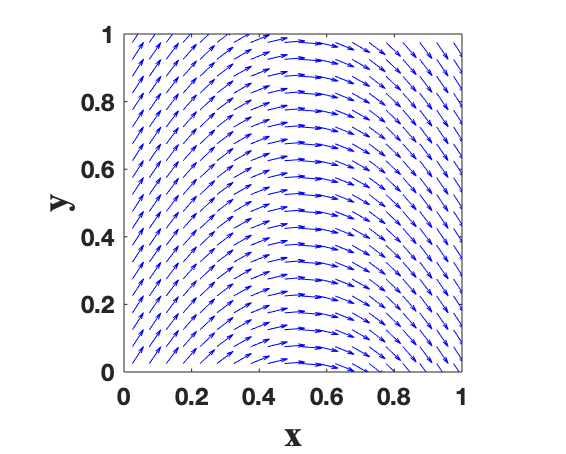}}
    \subfloat[angle profile]{\includegraphics[width=0.53\linewidth]{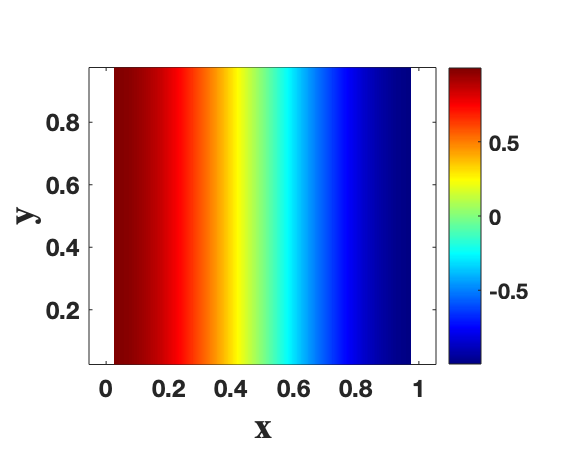}}
    \hspace{0.1in}
    \subfloat[arrow profile]{\includegraphics[width=0.5\linewidth]{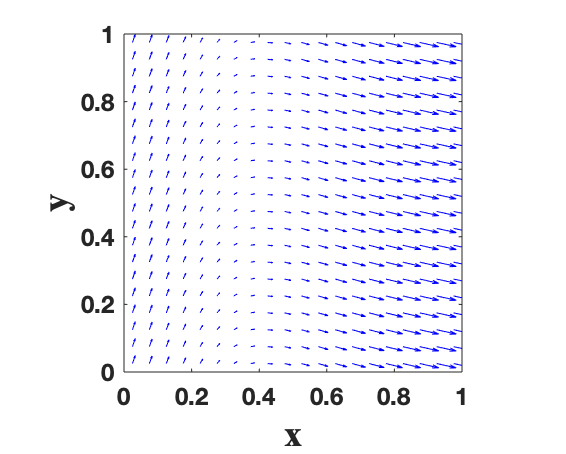}}
    \subfloat[angle profile]{\includegraphics[width=0.53\linewidth]{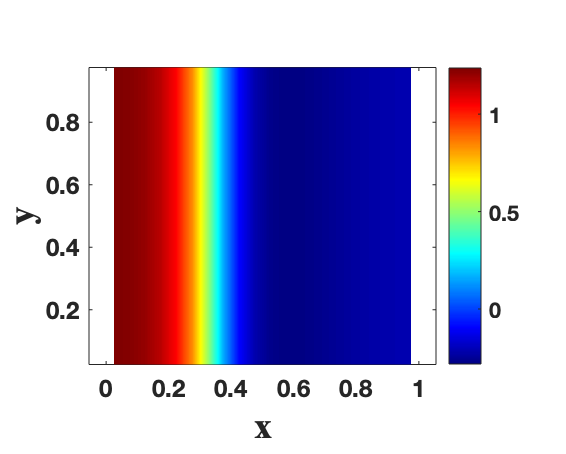}}
    \hspace{0.1in}
    \subfloat[arrow profile]{\includegraphics[width=0.5\linewidth]{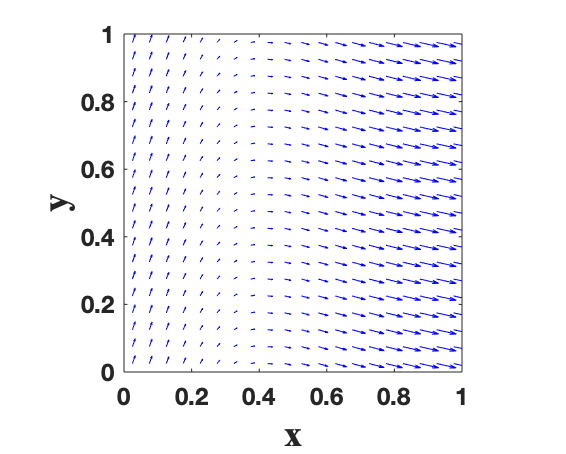}}
    \subfloat[angle profile]{\includegraphics[width=0.53\linewidth]{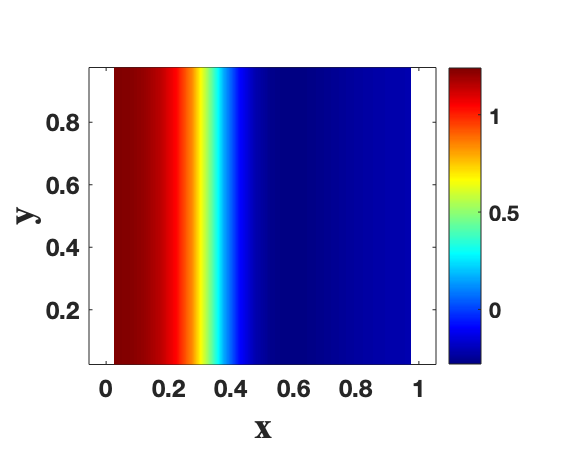}}
      \caption{The solution profile using GSPM and proposed methods in 3D given the initial condition $m_0$ with initial condition specified without source term, $\alpha=0$ and $T=0.1$, $N_x=N_y=N_z=20$, $N_t=400$. Top row with initial condition; Middle row with GSPM; Bottom row with proposed method. Initial condition given: $\m_0=[\cos(\cos(\pi x))\sin(0.01),\sin(\cos(\pi x))\sin(0.01),\cos(0.01)]$}
    \label{fig:3}
\end{figure}

\begin{figure}[htbp]
    \centering
    \subfloat[arrow profile]{\includegraphics[width=0.5\linewidth]{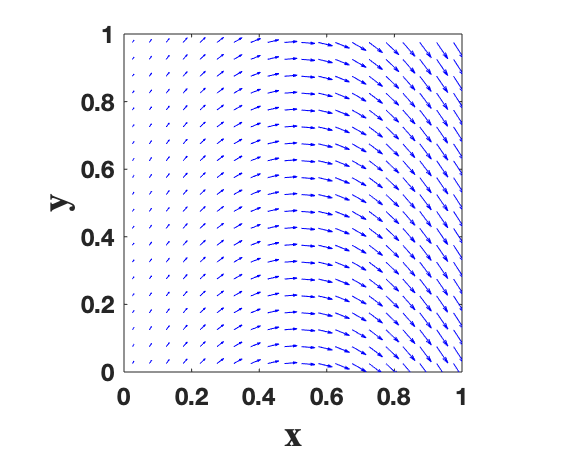}}
    \subfloat[angle profile]{\includegraphics[width=0.53\linewidth]{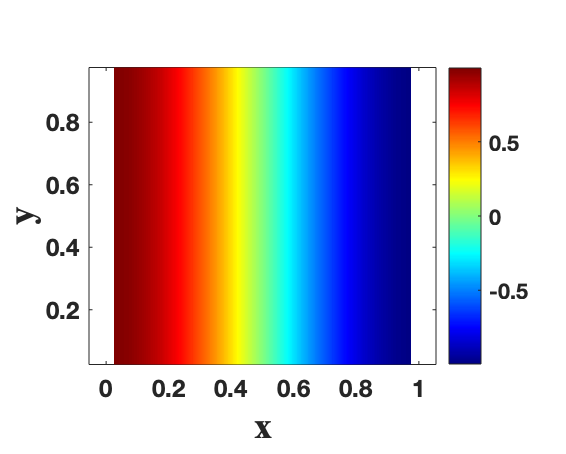}}
    \hspace{0.1in}
    \subfloat[arrow profile]{\includegraphics[width=0.5\linewidth]{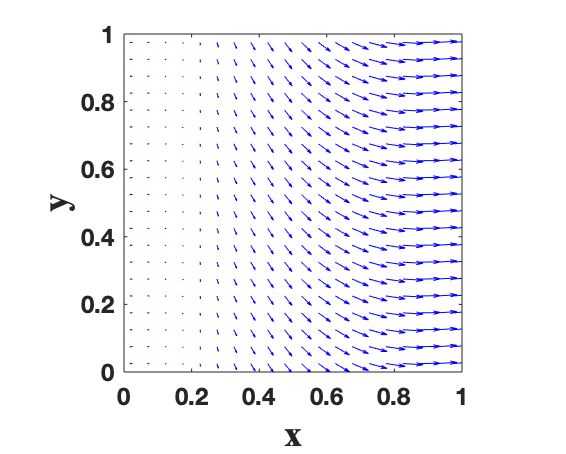}}
    \subfloat[angle profile]{\includegraphics[width=0.53\linewidth]{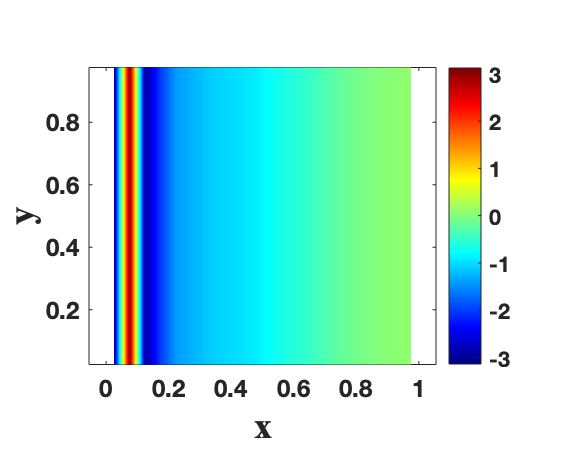}}
    \hspace{0.1in}
    \subfloat[arrow profile]{\includegraphics[width=0.5\linewidth]{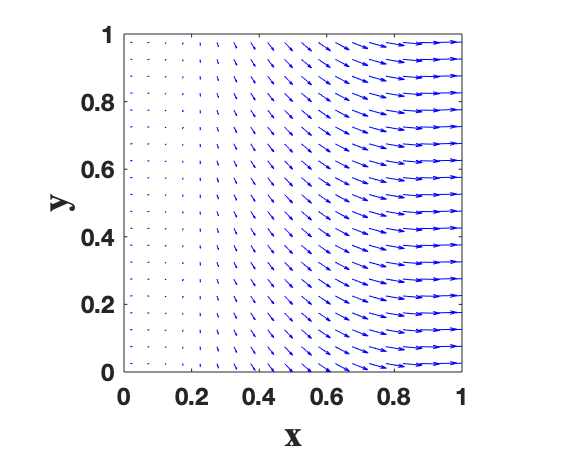}}
    \subfloat[angle profile]{\includegraphics[width=0.53\linewidth]{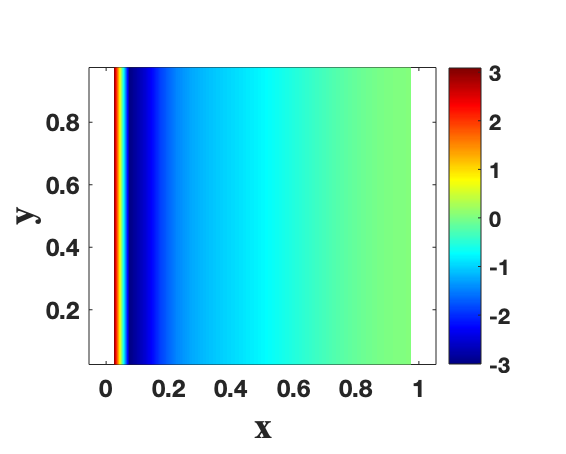}}
      \caption{The solution profile using GSPM and proposed methods in 3D given the initial condition $m_0$ with initial condition specified without source term, $\alpha=0$ and $T=0.1$, $N_x=N_y=N_z=20$, $N_t=40$. Top row with initial condition; Middle row with GSPM; Bottom row with proposed method. Initial condition given: $\m_0=[\cos(\cos(\pi x))\sin(x+t),\sin(\cos(\pi x))\sin(x+t),\cos(x+t)]$ with $t=T0=0$.}
    \label{fig:4}
\end{figure}

\begin{figure}[htbp]
    \centering
    \subfloat[arrow profile]{\includegraphics[width=0.5\linewidth]{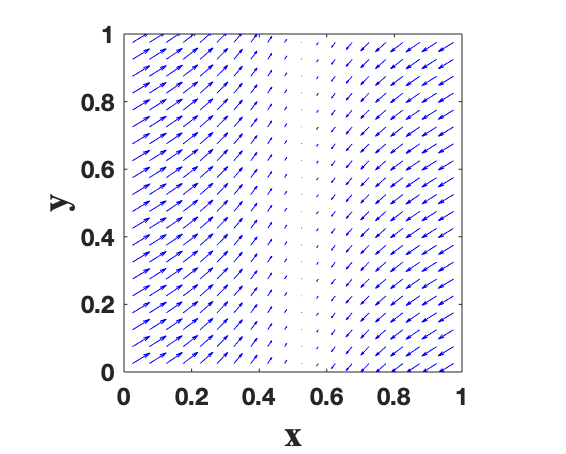}}
    \subfloat[angle profile]{\includegraphics[width=0.53\linewidth]{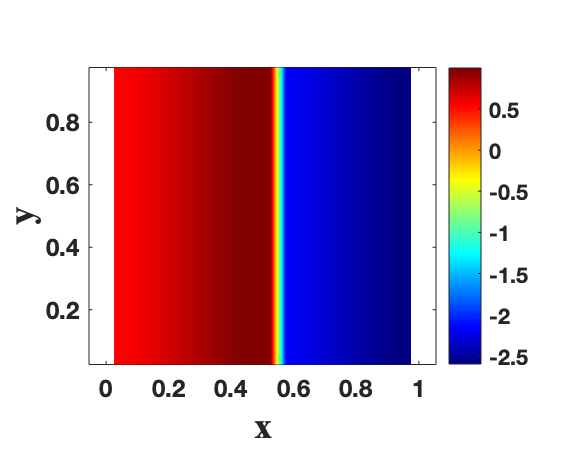}}
    \hspace{0.1in}
    \subfloat[arrow profile]{\includegraphics[width=0.5\linewidth]{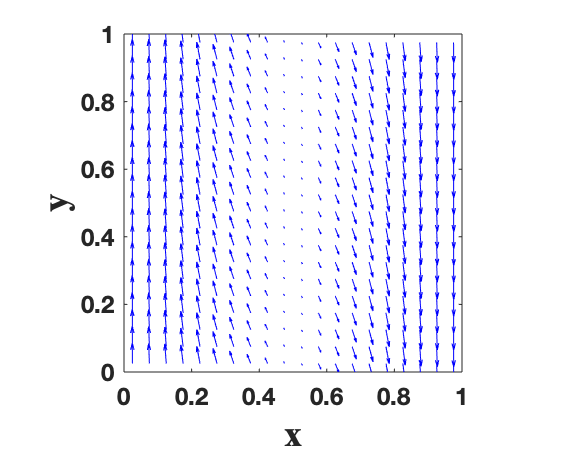}}
    \subfloat[angle profile]{\includegraphics[width=0.53\linewidth]{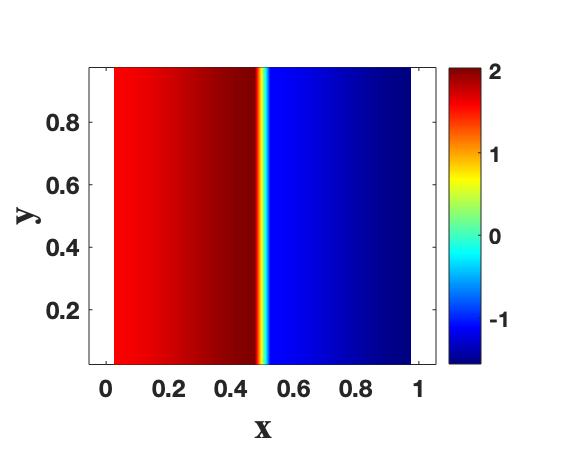}}
    \hspace{0.1in}
    \subfloat[arrow profile]{\includegraphics[width=0.5\linewidth]{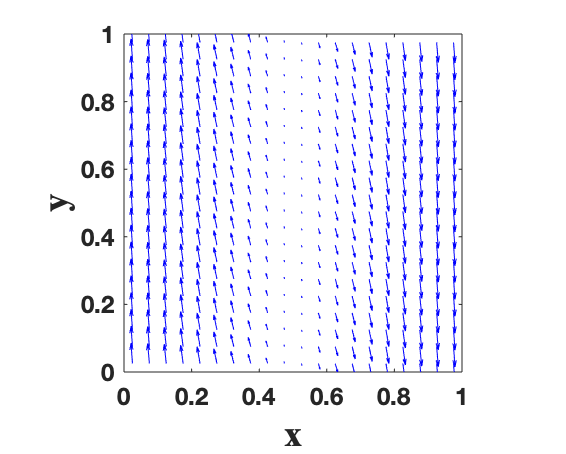}}
    \subfloat[angle profile]{\includegraphics[width=0.53\linewidth]{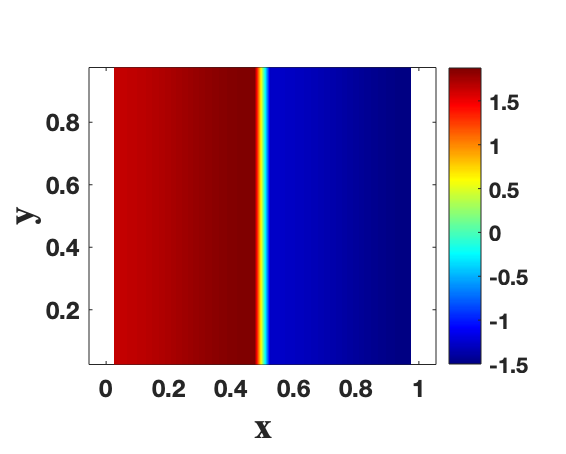}}
      \caption{The solution profile using GSPM and proposed methods in 3D given the initial condition $m_0$ with initial condition specified without source term, $\alpha=0$ and $T=0.1$, $N_x=N_y=N_z=20$, $N_t=40$. Top row with initial condition; Middle row with GSPM; Bottom row with proposed method. Initial condition given: $\m_0=[\cos(\cos(\cos(\pi x)))\sin(\pi x+t),\sin(\cos(\cos(\pi x)))\sin(\pi x+t),\cos(\pi x+t)]$ with $t=T0=0$.}
    \label{fig:5}
\end{figure}

\begin{figure}[htbp]
    \centering
    \subfloat[arrow profile]{\includegraphics[width=0.5\linewidth]{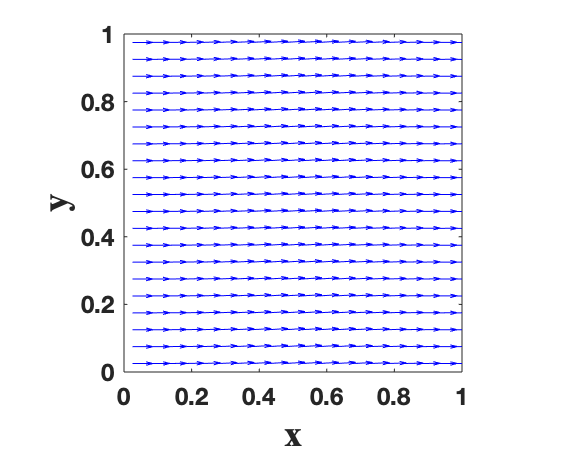}}
    \subfloat[angle profile]{\includegraphics[width=0.53\linewidth]{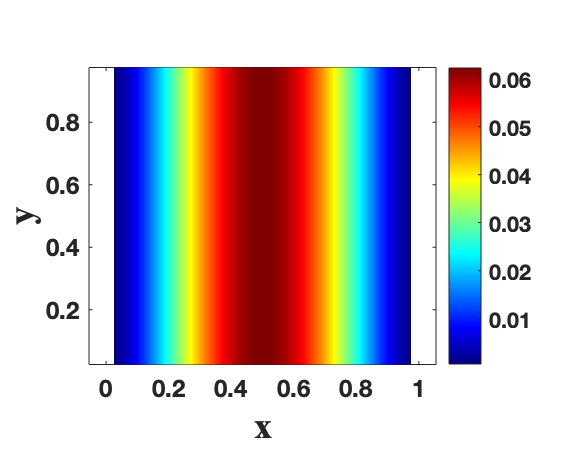}}
    \hspace{0.1in}
    \subfloat[arrow profile]{\includegraphics[width=0.5\linewidth]{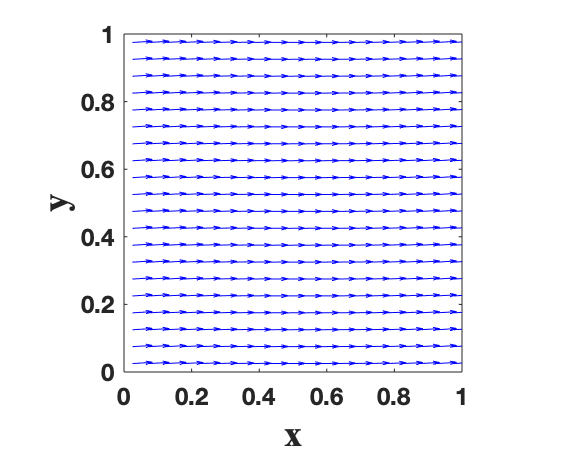}}
    \subfloat[angle profile]{\includegraphics[width=0.53\linewidth]{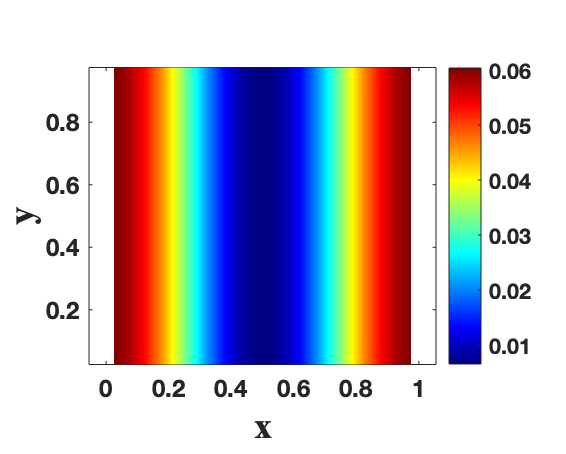}}
    \hspace{0.1in}
    \subfloat[arrow profile]{\includegraphics[width=0.5\linewidth]{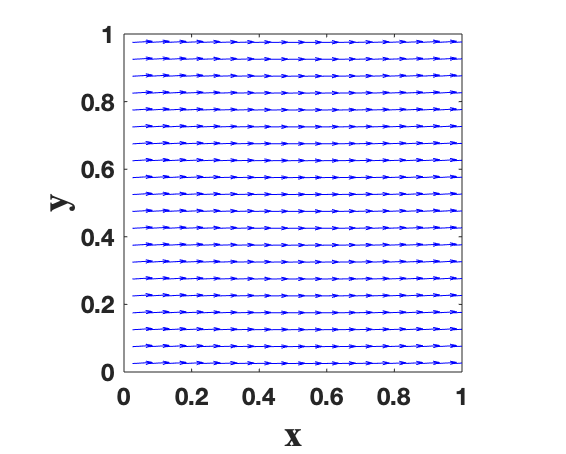}}
    \subfloat[angle profile]{\includegraphics[width=0.53\linewidth]{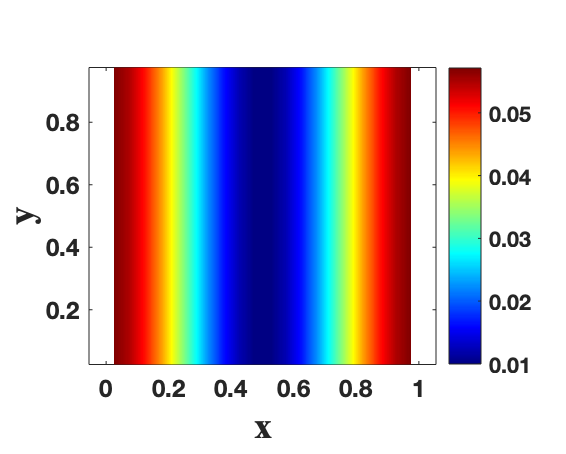}}
      \caption{The solution profile using GSPM and proposed methods in 3D given the initial condition $m_0$ with initial condition specified without source term, $\alpha=0$ and $T=0.1$, $N_x=N_y=N_z=20$, $N_t=40$. Top row with initial condition; Middle row with GSPM; Bottom row with proposed method. Initial condition given: $\m_0=[\cos(X)\sin(0.01),\sin(X)\sin(0.01),\cos(0.01)]$}
    \label{fig:6}
\end{figure}

In 3D, we choose the initial conditions below,
\begin{align*}
    \m_0&=\left(\cos(\cos(\pi x))\sin (0.01), \sin(\cos(\pi x))\sin (0.01), \cos (0.01)\right)^T\\
    \m_0&=\left(\cos(\cos(\pi x))\sin (x+t), \sin(\cos(\pi x))\sin (x+t), \cos (x+t)\right)^T\\
    \m_0&=\left(\cos(\cos(\cos(\pi x)))\sin (\pi x+t), \sin(\cos(\cos(\pi x)))\sin (\pi x+t), \cos (\pi x+t)\right)^T\\
     \m_0&=\left(\cos(X)\sin (0.01), \sin(X)\sin (0.01), \cos (0.01)\right)^T.
\end{align*}
The results for those initial conditions are presented in \Cref{fig:3}, \Cref{fig:4}, \Cref{fig:5}, and \Cref{fig:6} which verify the consistency.

\begin{figure}[htbp]
    \centering
    \subfloat[arrow profile]{\includegraphics[width=0.5\linewidth]{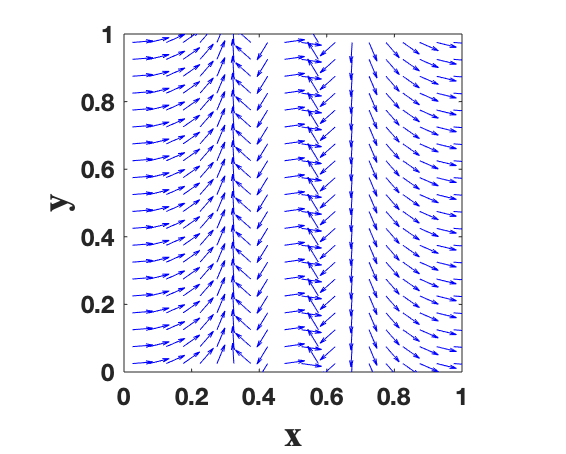}}
    \subfloat[angle profile]{\includegraphics[width=0.53\linewidth]{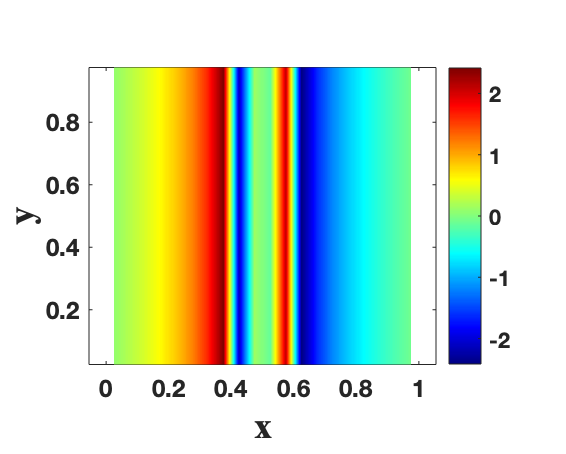}}
    \hspace{0.1in}
    \subfloat[arrow profile]{\includegraphics[width=0.5\linewidth]{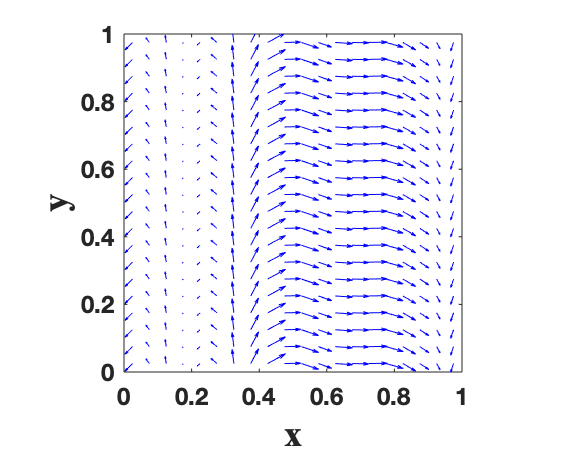}}
    \subfloat[angle profile]{\includegraphics[width=0.53\linewidth]{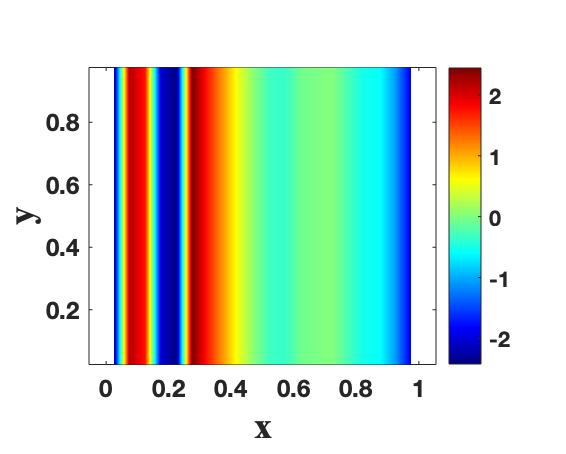}}
    \hspace{0.1in}
    \subfloat[arrow profile]{\includegraphics[width=0.5\linewidth]{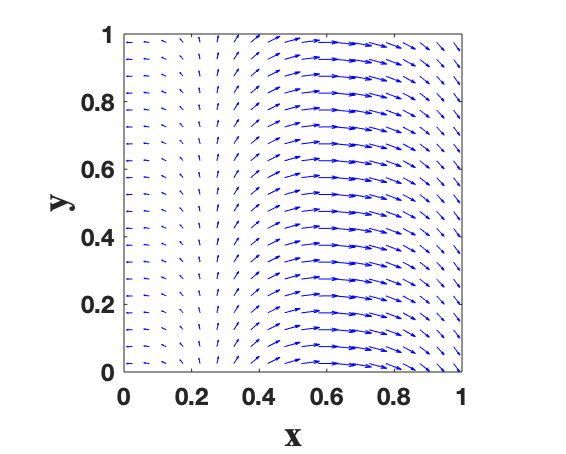}}
    \subfloat[angle profile]{\includegraphics[width=0.53\linewidth]{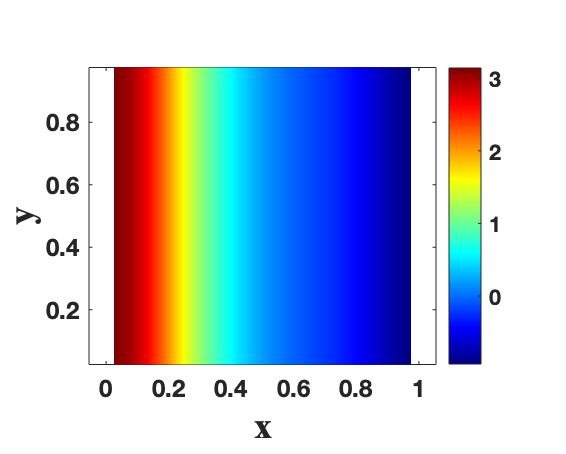}}
      \caption{The solution profile using GSPM and proposed methods in 3D given the initial condition $m_0$ with initial condition specified without source term, $\alpha=0$ and $T=0.1$, $N_x=N_y=N_z=20$, $N_t=40$. Top row with initial condition; Middle row with GSPM; Bottom row with proposed method. Initial condition given: $\m_0=[\cos(\tan(\pi x))\sin(0.01),\sin(\tan(\pi x))\sin(0.01),\cos(0.01)]$}
    \label{fig:7}
\end{figure}

\begin{figure}[htbp]
    \centering
    \subfloat[arrow profile]{\includegraphics[width=0.5\linewidth]{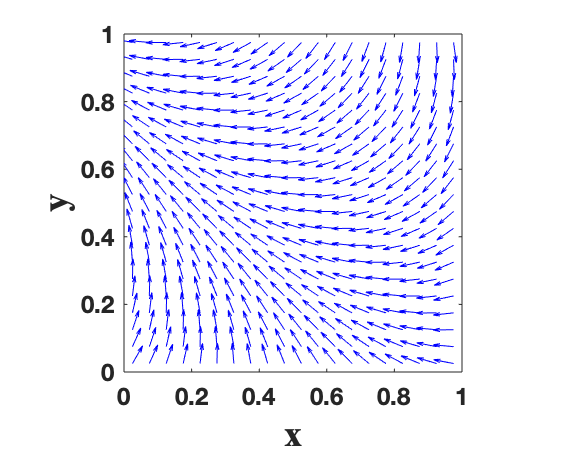}}
    \subfloat[angle profile]{\includegraphics[width=0.53\linewidth]{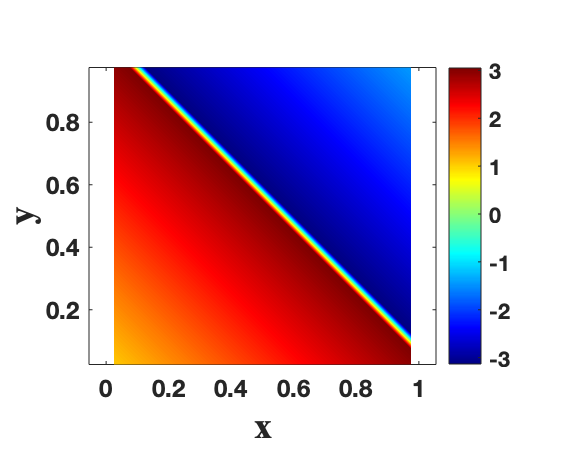}}
    \hspace{0.1in}
    \subfloat[arrow profile]{\includegraphics[width=0.5\linewidth]{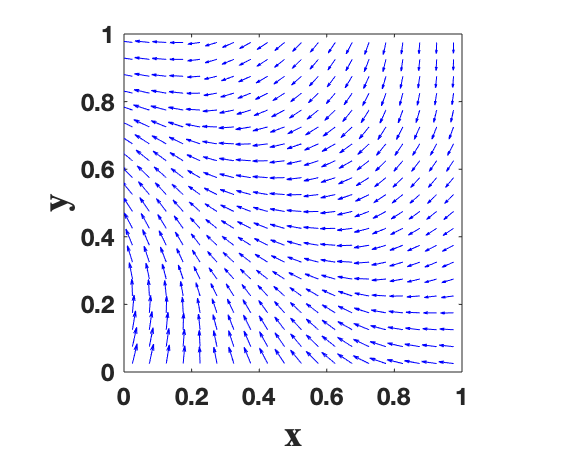}}
    \subfloat[angle profile]{\includegraphics[width=0.53\linewidth]{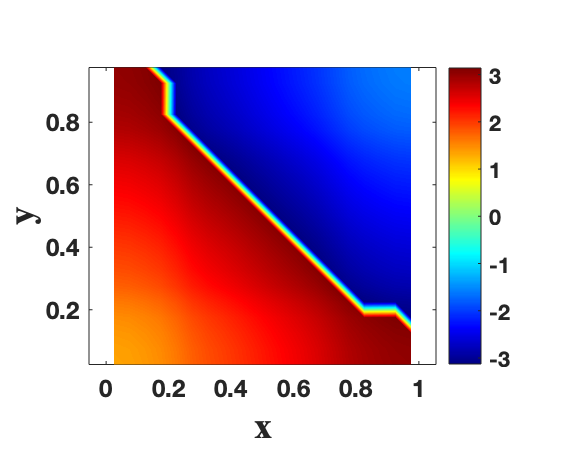}}
    \hspace{0.1in}
    \subfloat[arrow profile]{\includegraphics[width=0.5\linewidth]{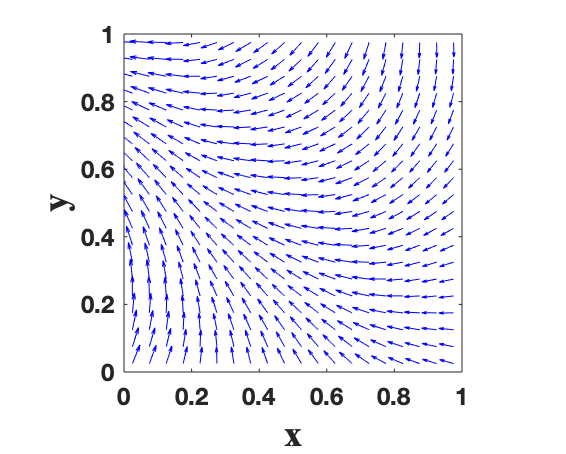}}
    \subfloat[angle profile]{\includegraphics[width=0.53\linewidth]{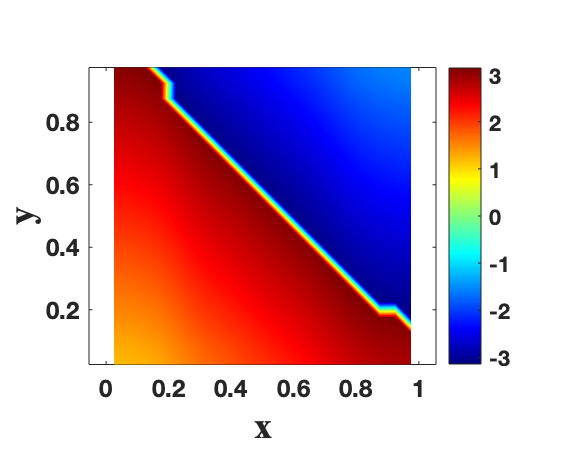}}
      \caption{The solution profile using GSPM and proposed methods in 3D given the initial condition $m_0$ with initial condition specified without source term, $\alpha=0$ and $T=0.01$, $N_x=N_y=N_z=20$, $N_t=40$. Top row with initial condition; Middle row with GSPM; Bottom row with proposed method. Initial condition given: $\m_0=[\cos(2(x+y+z))\sin(0.01),\sin(2(x+y+z))\sin(0.01),\cos(0.01)]$}
    \label{fig:8}
\end{figure}

\begin{figure}[htbp]
    \centering
    \subfloat[arrow profile]{\includegraphics[width=0.5\linewidth]{v1/init_3_angle_v1.png}}
    \subfloat[angle profile]{\includegraphics[width=0.53\linewidth]{v1/init_3_color_v1.png}}
    \hspace{0.1in}
    \subfloat[arrow profile]{\includegraphics[width=0.5\linewidth]{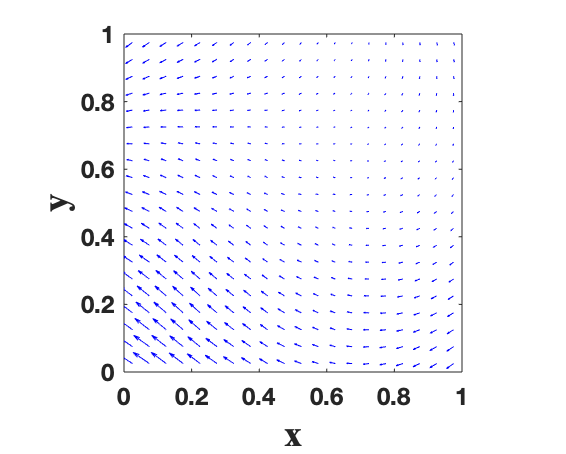}}
    \subfloat[angle profile]{\includegraphics[width=0.53\linewidth]{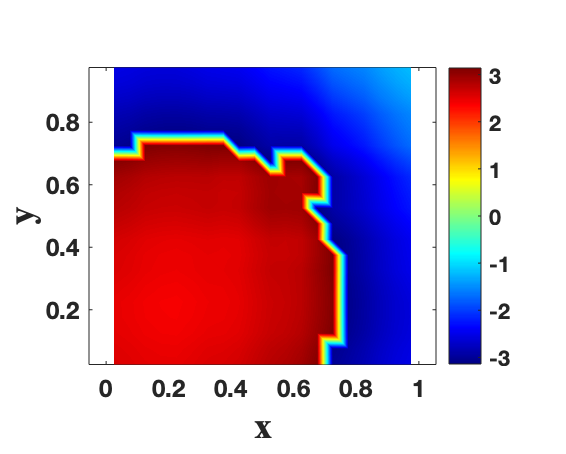}}
    \hspace{0.1in}
    \subfloat[arrow profile]{\includegraphics[width=0.5\linewidth]{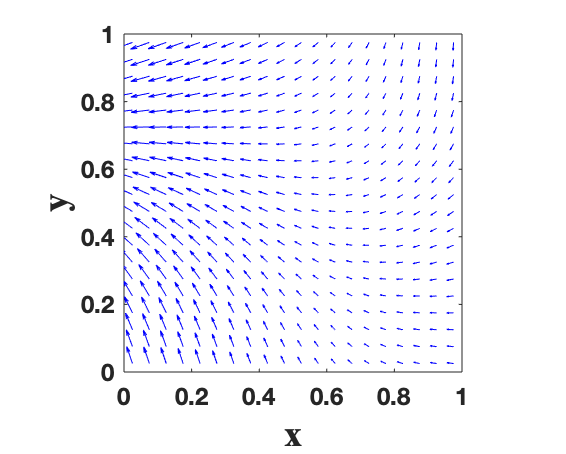}}
    \subfloat[angle profile]{\includegraphics[width=0.53\linewidth]{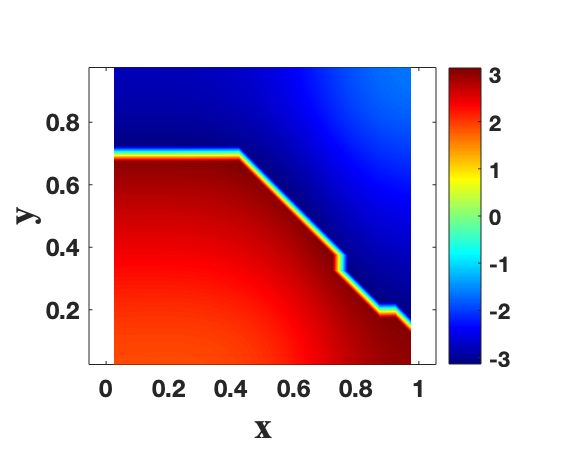}}
      \caption{The solution profile using GSPM and proposed methods in 3D given the initial condition $m_0$ with initial condition specified without source term, $\alpha=0$ and $T=0.1$, $N_x=N_y=N_z=20$, $N_t=40$. Top row with initial condition; Middle row with GSPM; Bottom row with proposed method. Initial condition given: $\m_0=[\cos(2(x+y+z))\sin(0.01),\sin(2(x+y+z))\sin(0.01),\cos(0.01)]$}
    \label{fig:9}
\end{figure}

\begin{figure}[htbp]
    \centering
    \subfloat[arrow profile]{\includegraphics[width=0.5\linewidth]{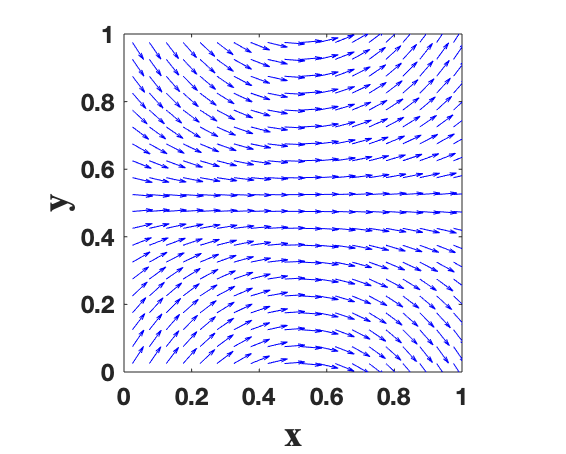}}
    \subfloat[angle profile]{\includegraphics[width=0.53\linewidth]{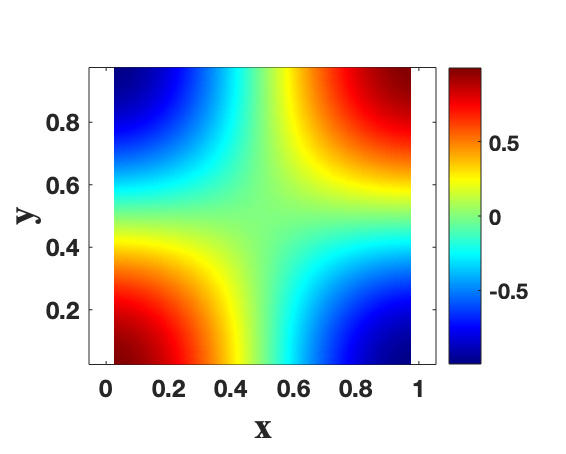}}
    \hspace{0.1in}
    \subfloat[arrow profile]{\includegraphics[width=0.5\linewidth]{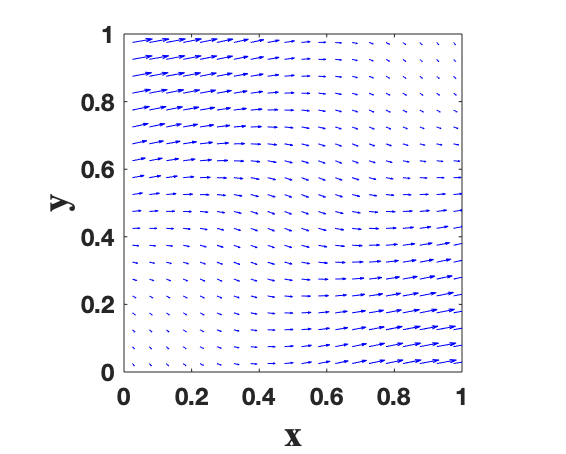}}
    \subfloat[angle profile]{\includegraphics[width=0.53\linewidth]{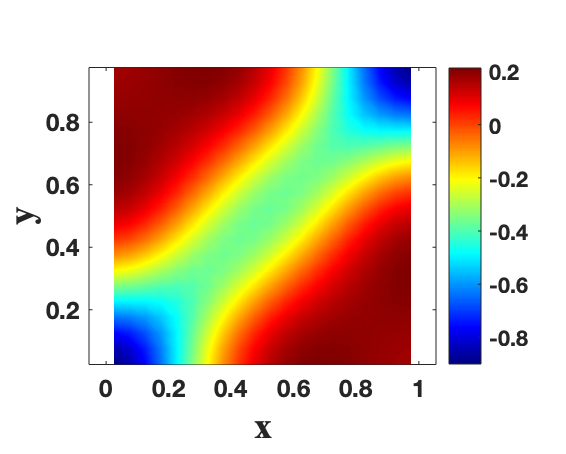}}
    \hspace{0.1in}
    \subfloat[arrow profile]{\includegraphics[width=0.5\linewidth]{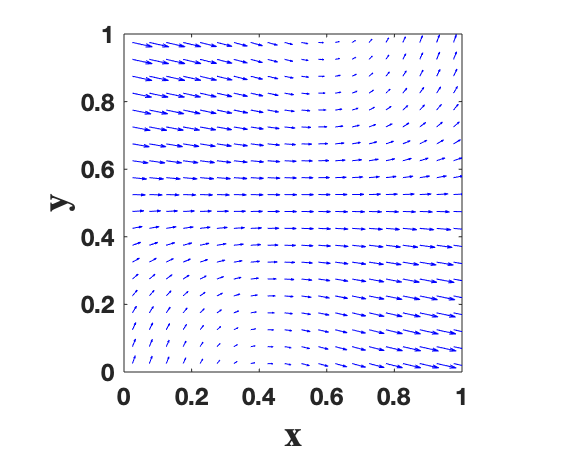}}
    \subfloat[angle profile]{\includegraphics[width=0.53\linewidth]{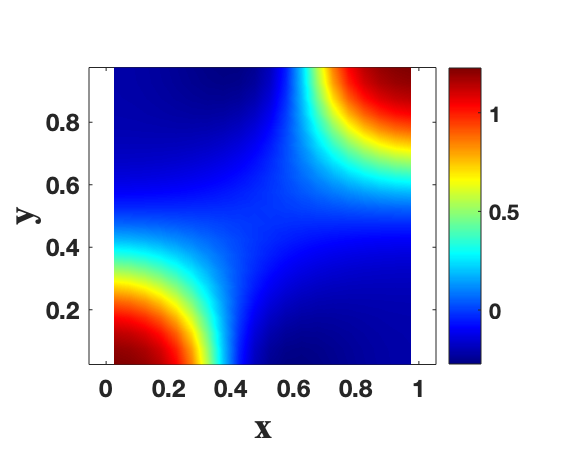}}
      \caption{The solution profile using GSPM and proposed methods in 3D given the initial condition $m_0$ with initial condition specified without source term, $\alpha=0$ and $T=0.1$, $N_x=N_y=N_z=20$, $N_t=400$. Top row with initial condition; Middle row with GSPM; Bottom row with proposed method. Initial condition given: $\m_0=[\cos(\cos(\pi x)\cos(\pi y))\sin(0.01),\sin(\cos(\pi x)\cos(\pi y))\sin(0.01),\cos(0.01)]$}
    \label{fig:10}
\end{figure}

\begin{figure}[htbp]
    \centering
    \subfloat[arrow profile]{\includegraphics[width=0.5\linewidth]{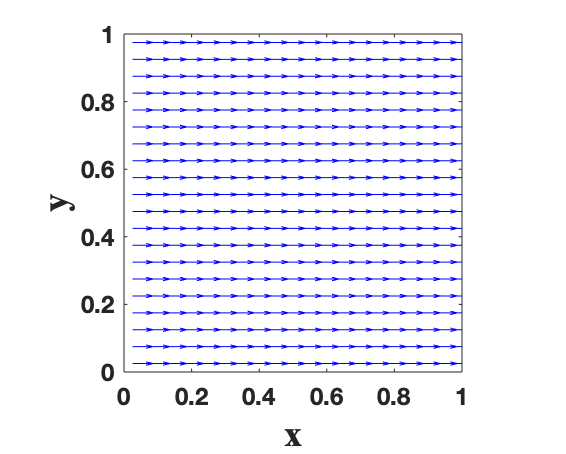}}
    \subfloat[angle profile]{\includegraphics[width=0.53\linewidth]{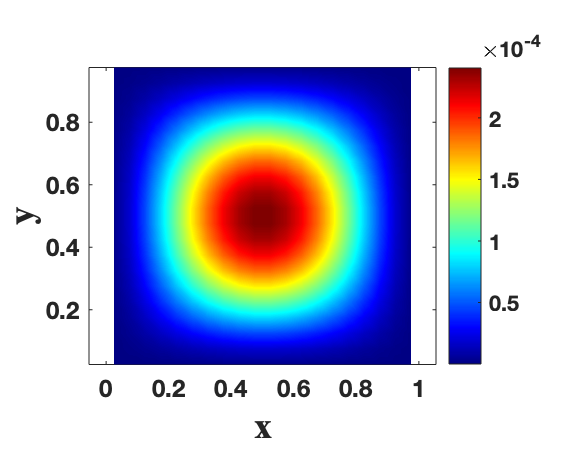}}
    \hspace{0.1in}
    % \subfloat[arrow profile]{\includegraphics[width=0.5\linewidth]{v1/initial_3D_arrow_v2.png}}
    % \subfloat[angle profile]{\includegraphics[width=0.53\linewidth]{v1/initial_3D_color_v2.png}}
    \subfloat[arrow profile]{\includegraphics[width=0.5\linewidth]{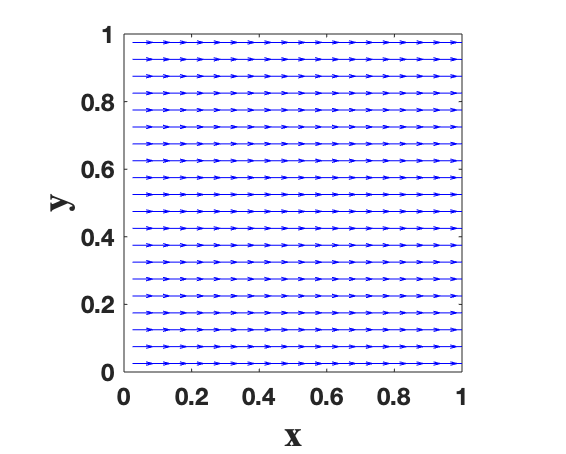}}
    \subfloat[angle profile]{\includegraphics[width=0.53\linewidth]{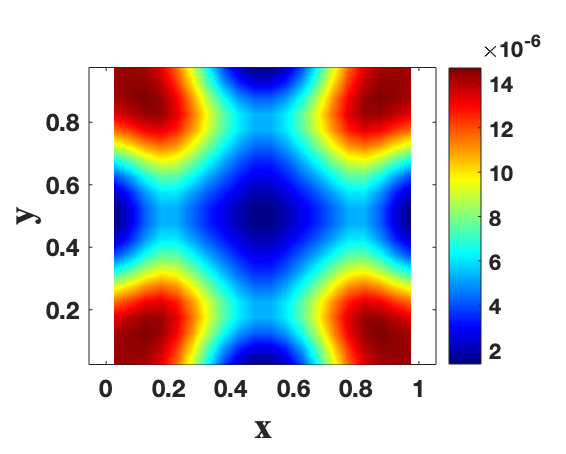}}
    \hspace{0.1in}
    \subfloat[arrow profile]{\includegraphics[width=0.5\linewidth]{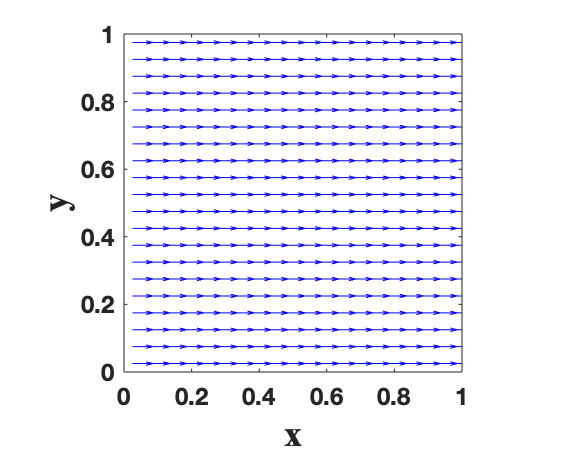}}
    \subfloat[angle profile]{\includegraphics[width=0.53\linewidth]{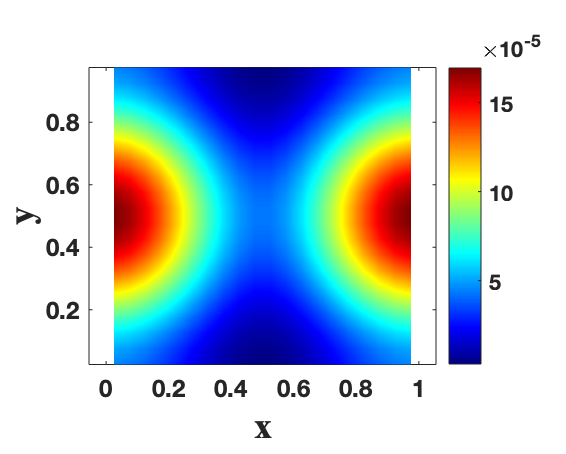}}
      \caption{Given initial conditions in 3D as $\m_0=[\cos(XYZ)\sin(0.01),\sin(XYZ)\sin(0.01),\cos(0.01)]$ specified without source term, $\alpha=0.01$ and $T=0.1$, $N_x=N_y=N_z=20$, $N_t=40$. Top row: initial condition; Middle row: GSPM; Bottom row: proposed method. }
    \label{fig:11}
\end{figure}

To get the robustness, we test other initial conditions below,
\begin{align*}
    \m_0&=\left(\cos(\tan(\pi x))\sin (0.01), \sin(\tan(\pi x))\sin (0.01), \cos (0.01)\right)^T\\
     \m_0&=\left(\cos(2(x+y+z))\sin (0.01), \sin(2(x+y+z))\sin (0.01), \cos (0.01)\right)^T\\
       \m_0&=\left(\cos(\cos(\pi x)\cos(\pi y))\sin (0.01), \sin(\cos(\pi x)\cos(\pi y))\sin (0.01), \cos (0.01)\right)^T\\
        \m_0&=\left(\cos(XYZ)\sin (0.01), \sin(XYZ)\sin (0.01), \cos (0.01)\right)^T\\
         \m_0&=\left(\cos(\cos(\pi x)\cos(\pi y)\cos(\pi z))\sin (0.01), \sin(\cos(\pi x)\cos(\pi y)\cos(\pi z))\sin (0.01), \cos (0.01)\right)^T.
\end{align*}
The results for other initial conditions are presented in \Cref{fig:7}, \Cref{fig:8}, \Cref{fig:9}, \Cref{fig:10} and \Cref{fig:11} which verify the stability and robustness.

\begin{figure}[htbp]
    \centering
    % \subfloat[arrow profile]{\includegraphics[width=0.5\linewidth]{v1/initial_3D_arrow_v1.png}}
    % \subfloat[angle profile]{\includegraphics[width=0.53\linewidth]{v1/initial_3D_color_v1.png}}
    % \hspace{0.1in}
    \subfloat[arrow profile]{\includegraphics[width=0.5\linewidth]{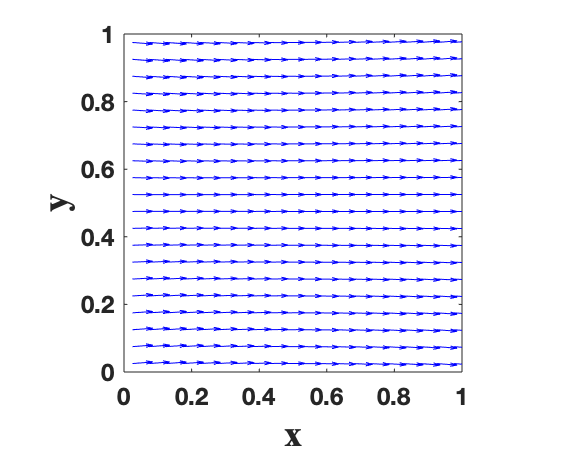}}
    \subfloat[angle profile]{\includegraphics[width=0.53\linewidth]{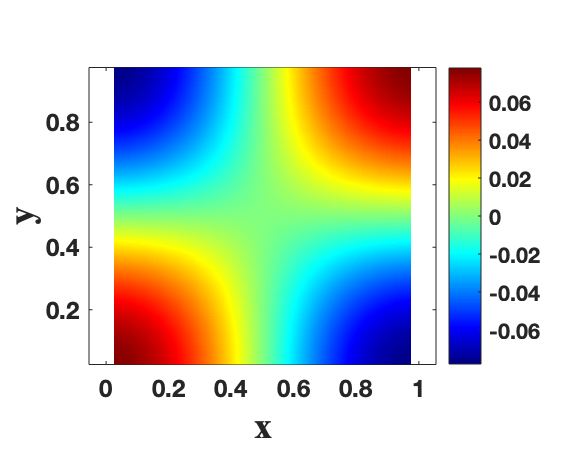}}
     \hspace{0.1in}
    \subfloat[arrow profile]{\includegraphics[width=0.5\linewidth]{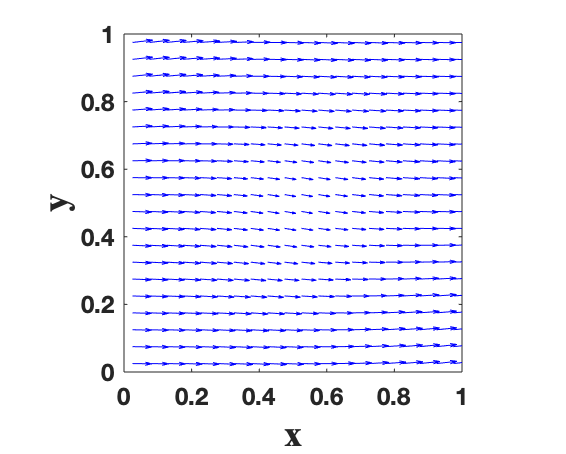}}
    \subfloat[angle profile]{\includegraphics[width=0.53\linewidth]{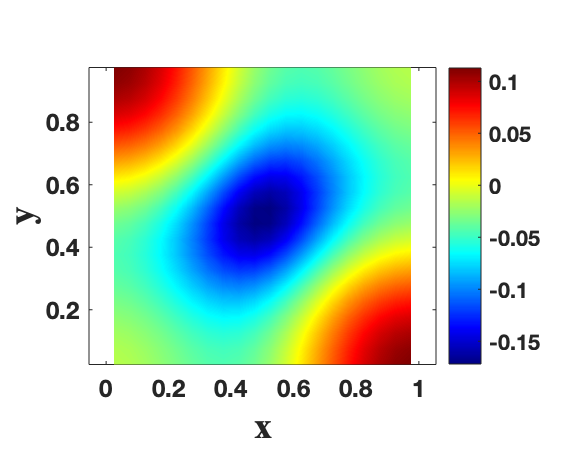}}
    \hspace{0.1in}
    \subfloat[arrow profile]{\includegraphics[width=0.5\linewidth]{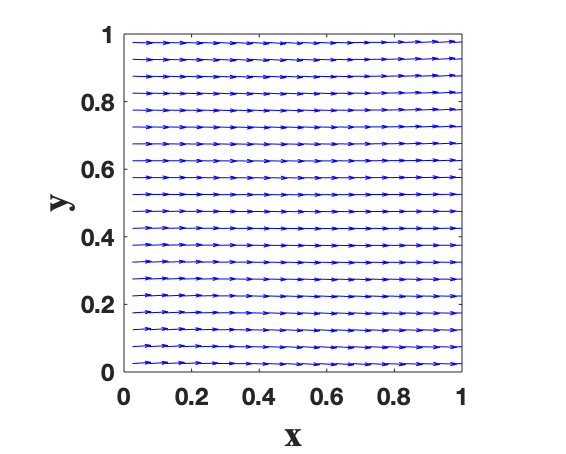}}
    \subfloat[angle profile]{\includegraphics[width=0.53\linewidth]{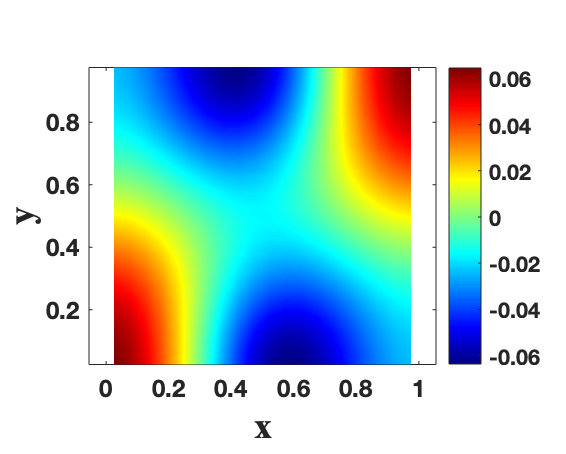}}
      \caption{Given initial conditions in 3D as $\m_0=[\cos(\cos(\pi x)\cos(\pi y)\cos(\pi z))\sin(0.01),\sin(\cos(\pi x)\cos(\pi y)\cos(\pi z))\sin(0.01),\cos(0.01)]$ specified without source term, $\alpha=0.01$ and $T=0.1$, $N_x=N_y=N_z=20$, $N_t=40$. Top row: initial condition; Middle row: GSPM; Bottom row: proposed method.}
    \label{fig:12}
\end{figure}

% \subsection{Stability tests}

% \subsection{Magnetic switching}

\section{Conclusions and discussions}
\label{sec:conclusions}

In this paper, we propose a structure preserving method with first order accuracy in time and second order accuracy in space. Such a method combines an implicit Gauss-Seidel iteration, a double diffusion iteration and a Crank-Nicolson's type iteration. The frist two step with an implicit Gauss-Seidel iteration, a double diffusion iteration gives a stable configuration for the higher order term. The last step to preserve the norm constraint. Such a method is constructed only based on the equation, not using a projection step. In the future work, we will analyze the stability of such a method and apply for micromagnetics simulations.

\section*{Acknowledgments}
This work is supported in part by the Jiangsu Science and Technology Programme-Fundamental Research Plan Fund (BK20250468), Research and
Development Fund of Xi'an Jiaotong Liverpool University (RDF-24-01-015).

\vspace{1cm}

\bibliographystyle{elsarticle-num-names}
\bibliography{references.bib}

\end{document}